\documentclass[12pt,reqno]{amsart}
\usepackage{amssymb,amscd,amsthm}
\usepackage[T1]{fontenc}
\usepackage{blkarray}

\usepackage{amsmath,amssymb,graphicx,mathrsfs} 
\usepackage[colorlinks=true,urlcolor=blue,citecolor=red,linkcolor=blue]{hyperref}

\let\frak\mathfrak
\let\Bbb\mathbb

\def\>{\relax\ifmmode\mskip.666667\thinmuskip\relax\else\kern.111111em\fi}
\def\<{\relax\ifmmode\mskip-.333333\thinmuskip\relax\else\kern-.0555556em\fi}
\def\vsk#1>{\vskip#1\baselineskip}
\def\vv#1>{\vadjust{\vsk#1>}\ignorespaces}
\def\vvn#1>{\vadjust{\nobreak\vsk#1>\nobreak}\ignorespaces}

  \let\ssize\scriptstyle
\let\sssize\scriptscriptstyle

\let\Medskip\medskip
\def\medskip{\par\Medskip}
\let\Bigskip\bigskip
\def\bigskip{\par\Bigskip}

\let\Maketitle\maketitle
\def\maketitle{\Maketitle\thispagestyle{empty}\let\maketitle\empty}

\newtheorem{thm}{Theorem}[section]

\newtheorem{lem}[thm]{Lemma}
\newtheorem{prop}[thm]{Proposition}

\numberwithin{equation}{section}

\theoremstyle{definition}
\newtheorem{defn}[thm]{Definition}
\newtheorem{rem}[thm]{Remark}

\newtheorem*{example}{Example}

\let\mc\mathcal
\let\nc\newcommand

\let\al\alpha
\let\bt\beta
\let\dl\delta

\let\gm\gamma
\let\Gm\Gamma

\let\la\lambda
\let\La\Lambda

\let\phi\varphi

\let\sig\varsigma

\let\Ups\Upsilon
\let\om\omega

\let\der\partial

\let\geq\geqslant

\let\leq\leqslant

\let\on\operatorname
\let\bi\bibitem
\let\bs\boldsymbol

\def\C{{\mathbb C}}
\def\Z{{\mathbb Z}}
\def\Q{{\mathbb Q}}

\def\Pb{{\mathbb P}}

\def\F{{\mc F}}

\def\+#1{^{\{#1\}}}

\def\sln{\mathfrak{sl}_N}

\def\beq{\begin{equation}}
\def\eeq{\end{equation}}
\def\be{\begin{equation*}}
\def\ee{\end{equation*}}

\nc{\bea}{\begin{eqnarray*}}
\nc{\eea}{\end{eqnarray*}}
\nc{\bean}{\begin{eqnarray}}
\nc{\eean}{\end{eqnarray}}
\nc{\bal}{\begin{align*}}
\nc{\eal}{\end{align*}}
\nc{\baln}{\begin{align}}
\nc{\ealn}{\end{align}}

\nc{\Il}{{\mc I_{\bs\la}}}
\nc{\bla}{{\bs\la}}
\nc{\Fla}{\F_\bla}
\nc{\tfl}{{T^*\Fla}}
\nc{\GL}{{GL_n(\C)}}
\nc{\GLC}{{GL_n(\C)\times\C^*}}

\let\sd s 

\def\ddk_#1{\kk_{#1}\<\>\frac\der{\der\<\>\kk_{#1}}}

\def\bul{\mathbin{\raise.2ex\hbox{$\sssize\bullet$}}}
\def\intt{\mathchoice
{\mathop{\raise.2ex\rlap{$\,\,\ssize\backslash$}{\intop}}\nolimits}
{\mathop{\raise.3ex\rlap{$\,\sssize\backslash$}{\intop}}\nolimits}
{\mathop{\raise.1ex\rlap{$\sssize\>\backslash$}{\intop}}\nolimits}
{\mathop{\rlap{$\sssize\<\>\backslash$}{\intop}}\nolimits}}

\let\kk q 
\let\cc c

\let\Ko K

\def\GZ/{Gelfand-Zetlin}
\def\KZ/{{\slshape KZ\/}}
\def\qKZ/{{\slshape qKZ\/}}
\def\XXX/{{\slshape XXX\/}}

\nc{\slnl}{{\sln (\lambda)}}
\nc{\PCN}{{   (\C[x])^N   }}
\nc{\di}{\on{Diag}}
\nc{\dio}{\on{Diag}_0}
\nc{\Mm}{{\mc M}}
\nc{\Nn}{{\mc N}}

\nc{\A}{{\mc C}}

\nc{\PCr}{{  P  (\C[x])^n   }}

\nc{\Pk}{{(\bs{P}^1)^k}}

\nc{\N}{{\Bbb N}}

\nc{\Ll}{{\mc L}}

\nc{\ord}{{\on{ord}\,}}

\nc{\Sing}{{\on{Sing}\,}}
\nc{\sing}{{\on{Sing}\,}}

\nc{\Hess}{{\on{Hess}}}

\nc{\R}{{\Bbb R}}

\let\on\operatorname
\nc{\Kk}{{\bs K}}
\nc{\Ap}{{A_\Phi(z)}}
\nc{\ap}{{A_\Phi(z)}}

\nc{\sv}{{\sing V}}
\nc{\cd}{{\C^n-\Delta}}
\nc{\UT}{{U^0}}   
\nc{\ep}{\epsilon}

\usepackage[all,cmtip]{xy}
\usepackage{empheq}
\usepackage{bm}

\newlanguage\fakelanguage
\newcommand\cyr{\fontencoding{OT2}\fontfamily{wncyr}\selectfont
   \language\fakelanguage}
\DeclareTextFontCommand{\textcyr}{\cyr}

\numberwithin{equation}{section}

\usepackage{calligra}

\DeclareMathOperator{\HOM}{\mathscr{H}\text{\kern -3pt {\calligra\large om}}\,}

\makeatletter
\newsavebox{\@brx}
\newcommand{\llangle}[1][]{\savebox{\@brx}{\(\m@th{#1\langle}\)}%
  \mathopen{\copy\@brx\kern-0.5\wd\@brx\usebox{\@brx}}}
\newcommand{\rrangle}[1][]{\savebox{\@brx}{\(\m@th{#1\rangle}\)}%
  \mathclose{\copy\@brx\kern-0.5\wd\@brx\usebox{\@brx}}}
\makeatother

\usepackage{scalerel,stackengine}  
\stackMath
\newcommand\reallywidehat[1]{%
\savestack{\tmpbox}{\stretchto{%
  \scaleto{%
    \scalerel*[\widthof{\ensuremath{#1}}]{\kern.1pt\mathchar"0362\kern.1pt}%
    {\rule{0ex}{\textheight}}
  }{\textheight}%
}{2.4ex}}%
\stackon[-6.9pt]{#1}{\tmpbox}%
}
\newcommand{\bsh}{\begin{shaded}}
\newcommand{\esh}{\end{shaded}}

\newcommand{\DM}[2]{\overline{\mc M}_{#1,#2}}

\newcommand{\rqh}{|_{\substack{{\bf Q}=1\\ \hbar=1}}}
\newcommand{\edc}{\widehat\nabla}

\newcommand{\cem}{\gm_{\rm EM}}

\makeatletter
\@namedef{subjclassname@2020}{\textup{2020} Mathematics Subject Classification}
\makeatother

\begin{document}
\title[]{Borel $(\bm \al,\bm\bt)$-multitransforms and Quantum Leray--Hirsch: integral representations of solutions of quantum differential equations for $\Pb^1$-bundles}
\author[Giordano Cotti ]{Giordano Cotti$\>^\circ$ }
\keywords{Quantum cohomology, projective bundles, Borel integral transform}
\subjclass[2020]{53D45, 44A20, 44A30, 33B20}
\maketitle
\begin{center}
\textit{ $^\circ\>$Faculdade de Ci\^encias da Universidade de Lisboa\\ Grupo de F\'isica Matem\'atica \\
Campo Grande Edif\'icio C6, 1749-016 Lisboa, Portugal\/}

\end{center}
{\let\thefootnote\relax
\footnotetext{\vskip5pt 
\noindent
$^\circ\>$\textit{ E-mail}:  gcotti@fc.ul.pt, gcotti@sissa.it}}

\vskip6mm

\begin{abstract}
In this paper, we address the integration problem of the isomonodromic system of quantum differential equations ({\it qDE}s) associated with the quantum cohomology of $\mathbb P^1$-bundles on Fano varieties. It is shown that bases of solutions of the {\it qDE} of the total space of the $\mathbb P^1$-bundle can be reconstructed from the datum of bases of solutions of the corresponding {\it qDE} associated with the base space. This represents a quantum analog of the classical Leray--Hirsch theorem in the context of the isomonodromic approach to quantum cohomology.
The reconstruction procedure of the solutions can be performed in terms of some integral transforms, introduced in \cite{Cot22}, called {\it Borel $(\bm\al,\bm\bt)$-multitransforms}. We emphasize the emergence, in the explicit integral formulas, of an interesting sequence of special functions (closely related to iterated partial derivatives of the B\"ohmer--Tricomi incomplete Gamma function) as integral kernels. Remarkably, these integral kernels have a {\it universal} feature, being independent of the specifically chosen $\mathbb P^1$-bundle. When applied to projective bundles on products of projective spaces, our results give Mellin--Barnes integral representations of solutions of {\it qDE}s. As an example, we show how to integrate the {\it qDE} of blow-up of $\mathbb P^2$ at one point via Borel multitransforms of solutions of the {\it qDE} of $\mathbb P^1$.
\end{abstract}

\tableofcontents

\section{Introduction}

\noindent 1.1.\,\,Enumerative geometry is that branch of geometry that concerns the number of solutions to a given geometrical problem, rather than explicitly finding them all. In the last decades, ideas coming from physics brought innovation to enumerative geometry, with both new techniques and the emergence of new rich geometrical structures. As an example, Gromov--Witten theory, which focuses on counting numbers of curves on a target space, lead to the discovery of {\it quantum cohomology} and the closely related {\it quantum differential equations}. 

Given a smooth complex projective variety $X$, its quantum cohomology $QH^\bullet(X)$ is a family of commutative, associative, unital $\C$-algebra structures on $H^\bullet(X,\C)$, obtained by deforming the classical cohomological cup product. Such deformation is performed by adding some ``quantum correction terms'', containing information on the number of rational curves on $X$. Namely, the structure constants of the quantum cohomology algebras can be expressed as third derivatives of a generating power series $F^X_0(\bm t)$, with $\bm t=(t^1,\dots, t^n)$ and $n=\dim_\C H^\bullet(X,\C)$, of genus 0 Gromov--Witten invariants of $X$. Under the assumption of convergence of $F^X_0(\bm t)$ in some domain $M\subseteq H^\bullet(X,\C)\cong \C^n$, the points $\bm t\in M$ can be used to label the quantum algebra structures on $H^\bullet(X,\C)$, the corresponding product being denoted by $\circ_{\bm t}$. This equips the quantum cohomology $QH^\bullet(X)$ with an analytic Dubrovin--Frobenius structure, with $M$ being the underlying complex manifold \cite{Dub96, Man99, Her02, Sab08}.

Points $\bm t\in M$ are parameters of isomonodromic deformations of the {quantum differential equation} (for short, {\it qDE}) of $X$. This is a system of linear differential equations of the form
\beq\label{introeq1}
\frac{d}{{d}z}\sig(z, \bm t)=\left({\bm{\mc U}}(\bm t)+\frac{1}{z}\bm\mu(\bm t)\right)\sig(z,\bm t),
\eeq 
where $\sig$ is a $z$-dependent holomorphic vector field\footnote{Notice that tangent spaces of $M$ can canonically be identified with $H^\bullet(X,\C)$.} on $M$, and ${\bm{\mc U}},\bm \mu$ are two endomorphisms of the holomorphic tangent bundle of $M$. The first operator $\bm{\mc U}$ is the operator of $\circ$-multiplication by the Euler vector field, a distinguished vector field on $M$, obtained as perturbation of the constant vector field given by the first Chern class $c_1(X)$. The second operator $\bm\mu$, called {\it grading operator}, keeps track of the non-vanishing degrees of $H^\bullet(X,\C)$. 

The {\it qDE} is a rich object associated with $X$. In the first instance, the Gromov--Witten theory of $X$ can be reconstructed from the datum of the {\it qDE} \eqref{introeq1} only. For details on a Riemann--Hilbert--Birkhoff approach to reconstruct the generating function $F^X_0(\bm t)$, and consequently the Dubrovin--Frobenius structure of $QH^\bullet(X)$, see \cite{Dub96,Dub99,Cot21a,Cot21b}. In the second instance, the {\it qDE} of $X$ encodes not only information about the enumerative (or symplectic) geometry of $X$, but also (conjecturally) about its topology and complex geometry. In order to disclose such a great amount of information is via the study of the asymptotics and the monodromy of its solutions, see \cite{Dub98, GGI16, CDG18,Cot20}. The purpose of this paper is to construct new analytic tools, namely some integral representations of the solutions of {\it qDE}, which will be particularly convenient to the study of asymptotics, Stokes phenomena, and other analytical aspects. This will represent a continuation of the research direction started in \cite{Cot22}.
\vskip2mm

\noindent 1.2.\,\,The projectivization of vector bundles is one of the most natural constructions of smooth projective varieties. Projective bundles have consequently been among the first varieties whose quantum cohomology algebras have been studied. See \cite{QR98,CMR00,AM00}.

The role played by the Gromov--Witten theory of projective bundles has recently been revealed as central not only in the context of open deep conjectures (such as the crepant transformation conjecture for ordinary flops \cite{LLW16a,LLW16b,LLQW16}) but even for delicate foundational aspects of Gromov--Witten theory, such as its functoriality \cite{LLW15}.

The classical Leray--Hirsch theorem prescribes how to reconstruct the classical cohomology algebra $H^\bullet(P,\C)$ of a projective bundle $P=\Pb(V)\to X$ on a variety $X$, from the knowledge of the algebra $H^\bullet(X,\C)$ and the Chern roots $c_k(V)$, $k=0,\dots, {\rm rk }\,V$. These data only, indeed, allow to write down an explicit presentation of $H^\bullet(P,\C)$. Several ``quantum counterparts'' of this theorem have been proved over the years. Many of the main results proved in \cite{MP06,Ele05,Ele07,Bro14,LLW10,LLW15,LLW16a,LLW16b,Fan21} have a common thread: they allow to deduce information about the quantum cohomology (or, more generally, the Gromov--Witten theory) of a projective bundle $P\to X$ starting from information on the quantum cohomology of $X$. See Section \ref{secpb} for a more detailed discussion.

Following the same philosophy, in this paper we address the following question:

{\bf Q1.} {\it Is it possible to reconstruct a basis of solutions of the qDE of a projective bundle $P\to X$ from the datum of a basis of solutions of the qDE of $X$?}

We obtain a positive result, under the assumption that $P$ is a Fano split $\Pb^1$-bundle on $X$. In order to present the main results, we first briefly introduce some preliminary notions.
\vskip2mm

\noindent 1.3.\,\,The first notion we want to introduce is that of {\it master function}. We call master function of $X$ at $p\in QH^\bullet(X)$ any $\C$-valued function $\Phi_\sig$, holomorphic on the universal cover $\widetilde{\C^*}$ of $\C^*$, of the form
\[\Phi_\sig(z)=z^{-\frac{\dim_\C X}{2}}\int_X\sig(z,p),
\]where $\sig$ is a solution of the {\it qDE} \eqref{introeq1} specialized at $\bm t=p$.

Rather than working directly with the space of solutions of the {\it qDE}, it is more convenient to focus on the space $\mc S_p(X)$ of master functions of $X$ at $p\in QH^\bullet(X)$. More precisely, in addressing question {\bf Q1} above, we can split the problem into two parts:

{\bf Q2(1).} {\it Is it possible to reconstruct the space of solutions of the {\it qDE} (specialized at a point $p$) from the datum of the space $S_p(X)$ of master functions only?}

{\bf Q2(2).} {\it Is it possible to reconstruct the space of master functions $S_{\pi^*p}(P)$ from the datum of the space of master functions $S_p(X)$, where $p\in QH^\bullet(X)$ and\,\footnote{The map $\pi$ induces a map $\pi^*\colon H^\bullet(X,\C)\to H^\bullet(P,\C)$. The class $\pi^*p \in QH^\bullet(P)$ is the image of $p$.} $\pi\colon P\to X$?}

Question {\bf Q2(1)} has been extensively studied in \cite{Cot22}. In {\it loc.\,cit.}, it is shown that the answer is positive for generic $p\in QH^\bullet(X)$. The problem of reducing a vector differential equation to a scalar one is well known in the theory of ordinary differential equations. Such a scalar reduction is equivalent to the choice of what is traditionally called a {\it cyclic vector} for the differential system \cite[Lemma II.1.3]{Del70}. Moreover, several algorithmic reduction procedures have been developed, see e.g.\,\,\cite{Bar93} and references therein. On any Dubrovin--Frobenius manifold $M$, we have a natural candidate for the cyclic vector, namely the unit vector of the Frobenius algebras at each point $p\in M$. It turns out that such a choice is working on the complement of an analytic subset of $M$, called $\mc A_\La$-stratum \cite[Sec.\,2]{Cot22}. More details will be given in Section \ref{scyc}. Consequently, question {\bf Q2(2)} represents the main problem to be still addressed.

\begin{figure}[ht!]
\centering
\def\svgscale{1}
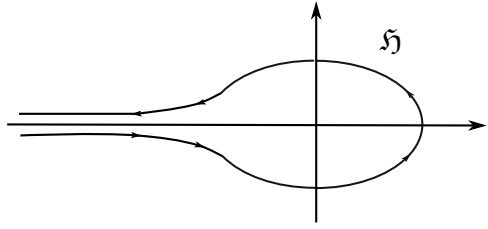
\caption{Hankel-type contour of integration defining Borel $(\bm \al,\bm\bt)$-multitransform.}
\label{gammahankel}
\end{figure}

The second notion we want to recall is that of ({\it analytic})  {\it Borel $(\bm\al,\bm\bt)$-multitransforms}, introduced in \cite{Cot22}. These are $\C$-multilinear integral transforms of tuples of analytic functions. Given two $h$-tuples $\bm \al,\bm\bt\in(\C^*)^{\times h}$, with $h\geq 1$, and an $h$-tuple of analytic functions $\Phi_h$ on $\C\setminus\R_{<0}$, we define its Borel $(\bm\al,\bm\bt)$-multitransform $\mathscr B_{\bm\al,\bm\bt}[\Phi_1,\dots,\Phi_h]$ by the integral
\[\mathscr B_{\bm\al,\bm\bt}\left[\Phi_1,\dots,\Phi_h\right](z):=\frac{1}{2\pi\sqrt{-1}}\int_{\frak H}\prod_{j=1}^h\Phi_j\left(z^\frac{1}{\al_j\bt_j}\la^{-\bt_j}\right)e^\la\frac{{\rm d}\la}{\la},
\]where $\frak H$ is a Hankel-type contour of integration, originating from $-\infty$, circling the origin once in the positive direction, and returning to $-\infty$. See Figure \ref{gammahankel}.

The third and last object we need to introduce for presenting our main results is a sequence of special functions $\mc E_k$, with $k\in\N$. Consider the function $\mc E(s,z)$, analytic and single valued on $\C\times{\widetilde{\C^*}}$, defined by the integral
\[\mc E(s,z):=\frac{z^s}{\Gm(s)}\int_0^1t^{s-1}e^{z-zt}{\rm d}t,\qquad{\rm Re}\,s>0, \quad z\in{\widetilde{\C^*}}.
\]Alternatively, the function $\mc E$ can be defined via the series expansion\footnote{Hence, $\mc E$ can be thought as a ``deformed'' exponential function. This explains the notation $\mc E$.}
\[\mc E(s,z)=\sum_{k=0}^\infty\frac{z^{k+s}}{\Gm(1+k+s)},\qquad (s,z)\in\C\times\widetilde{\C^*}.
\]
We define the functions $\mc E_k\in\mc O(\widetilde{\C^*})$, $k\in\N$, as the iterated partial derivatives
\[\mc E_k(z):=\left.\frac{\der^k}{\der s^k}\mc E(s,z)\right|_{s=0},\quad k\geq 0.
\]
For more explicit formulas for the $\mc E_k$-functions, see Section \ref{secek}.

The function $\mc E$ is closely related to the upper and lower incomplete Gamma functions. These are the functions $\Gm(s,z)$ and $\gm(s,z)$ defined by the integrals
\[\Gm(s,z):=\int_z^\infty t^{s-1}e^{-t}{\rm d}t,\qquad \gm(s,z):=\int_0^z t^{s-1}e^{-t}{\rm d}t,\qquad{\rm Re}\,s>0,\quad z\in\C.
\]The incomplete Gamma functions were first investigated,  for real $z$, in 1811 by A.M.\,Legendre \cite[Vol.\,1, pp.\,339--343 and later works]{Leg11}. The significance of the decomposition $\Gamma(s)=\gamma(s,z)+\Gamma(s,z)$ was recognised by F.A.\,Prym in 1877 \cite{Pry77}, who seems to have been the first to investigate the functional behavior of these functions (that he denoted by $P$ and $Q$, and which are often referred to as the {\it Prym functions}). The inconvenience of the function $\gamma(s,z)$ is not only of having poles at $s=0,-1,-2,\dots$, but even of being multivalued in $z$. Both inconveniences can be avoided by introducing, following P.E.\,B\"ohmer \cite{B\"oh39} and F.G.\,Tricomi \cite{Tri50}, the normalized incomplete Gamma function $\gamma^*$
\[\gamma^*(s,z):=\frac{z^{-s}}{\Gamma(s)}\gamma(s,z),\quad\text{which is entire on }\C^2.
\] Our function $\mc E$ satisfies the identities
\[\mc E(s,z)=e^z z^s \gamma^*(s,z)=e^z\frac{\gamma(s,z)}{\Gm(s)}=e^z\left(1-\frac{\Gm(s,z)}{\Gm(s)}\right).
\]
\begin{rem}
The classical theory of the incomplete Gamma functions (including series expansions of various kinds, asymptotics expansions, differentiation and recurrence relations, continued fractions, integral representations, etc.) can be found in \cite[Kap.\,II, XV,XXI]{Nie06} \cite[Kap.\,V]{B\"oh39}. This research was boostered after the late 1940s, when Tricomi recognised the importance of these functions\footnote{According to W.\,Gautschi, Tricomi was fascinated by the incomplete Gamma functions, and he was fond of calling them affectionately the {\it Cinderella functions}, see \cite{Gau98}.}, revitalized their study, and gave important contributions of his own. Tricomi himself summarized the knowledge as of 1950s in the second volume of the Bateman Project \cite[Ch.\,IX, pp.\,133--151]{Erd53}, and gave an even more detailed exposition in his monograph \cite[\S\S\,4.1--4.6]{Tri54}. For a beautiful overview of recent results on the incomplete Gamma functions since Tricomi see \cite{Gau98}. Curiously enough, however, in all these classical handbooks (and including \cite{AS64}, or the even more recent \cite{OLBC10}) the higher order derivatives $\der_s^k\Gm(s,z),\der_s^k\gm(s,z)$ and $\der^k_s\gm^* (s,z)$ are not studied for $k>1$. These exactly are those derivatives related to our functions $\mc E_k$. To the best of our knowledge, the first analysis of the functions $\der^k_s\Gm(s,z)$ was developed in \cite[\S 4]{GGMS90}.
\end{rem}
\vskip2mm
\noindent 1.4.\,\,In this paper, we give explicit integral representations of master functions of $P$ in terms of Borel $(\bm\al,\bm\bt)$-multitransforms of master functions of the base space $X$ under the following assumptions on $X$ and $P$:
\begin{enumerate}
\item[(I)] we assume that $X$ is a product $X=X_1\times\dots\times X_h$ of smooth Fano projective varieties $X_i$, 
\item[(II)] and that $\pi\colon P=\Pb(V)\to X$ is the projectivization of a split rank 2 vector bundle $V\cong \mc O_X\oplus L\to X$, where $L\to X$ is the external tensor product of fractional powers of the determinant bundles $\det TX_i$.
\end{enumerate}
Our first main result, Theorem \ref{mt2}, claims that any master function of $P$, at a point $\pi^*\bm\dl\in H^2(P,\C)$ of its small quantum cohomology, can be expressed in terms of Borel $(\bm\al,\bm\bt)$-multitransforms of master functions of $X_i$ at the point $\dl_i\in H^2(X_i,\C)$, where
$$\bm\dl=\sum_{j=1}^h1\otimes\dots\otimes\dl_j\otimes\dots\otimes 1.$$
More precisely, if $L=\boxtimes_{j=1}^hL_j^{\otimes(-d_j)}$ and $\det TX_j=L_j^{\otimes\ell_j}$, for ample line bundles $L_j\to X_j$ and $0<d_j<\ell_j$, then any master function of $P$ at $\pi^*\bm\dl$ is a $\C$-linear combination of integral of the form
\beq\label{borel1}
\mathscr B_{\bm\al,\bm\bt}\left[\Phi_1,\dots,\Phi_h,\mc E_k\right](z)=
\frac{1}{2\pi\sqrt{-1}}\int_{\frak H}\prod_{j=1}^h\Phi_j\left(z^{\frac{\ell_j-d_j}{\ell_j}}\la^\frac{d_j}{\ell_j}\right)\mc E_k\left(z^2\la^{-1}\right)e^\la\frac{{\rm d}\la}{\la},
\eeq
where
\begin{itemize}
\item $\bm\al,\bm\bt$ are the $(h+1)$-tuples
\[\bm\al=\left(\frac{\ell_1^2}{d_1(d_1-\ell_1)},\dots,\frac{\ell_h^2}{d_h(d_h-\ell_h)},\frac{1}{2}\right),
\qquad \bm\bt=\left(-\frac{d_1}{\ell_1},\dots,-\frac{d_h}{\ell_h},1\right),
\]
\item $\Phi_j$ is a master function of $X_j$ at the point $\dl_j\in H^2(X_j,\C)$,
\item and $k=0,\dots, \dim_\C X+1$.
\end{itemize}
\begin{rem}
We emphasize a {\it universal} feature of the integral kernels $\mc E_k$, their emergence in formula \eqref{borel1} being independent of the specifically chosen $\Pb^1$-bundle $P=\Pb(V)$ on $X$. 
\end{rem}

Our second main result concerns the more specific case of projective bundles on products of projective spaces. Consider the bundle $P=	\Pb(\mc O_X\oplus\mc O(-d_1,\dots, -d_h))$ on the variety $X=\Pb^{n_1-1}\times\dots\times\Pb^{n_h-1}$, with $0<d_i<n_i$ for $i=1,\dots, h$. Our Theorem \ref{mt2.2} claims that the space $\mc S_0(P)$ of master functions of $P$ at $0\in H^\bullet(P,\C)$ is a linear combination of integrals of the form
\[H_{\bm j}(z):=\frac{1}{(2\pi\sqrt{-1})^{h+1}}\int_{\times \gm_i}\int_\frak h\left[\prod_{i=1}^h\Gm(s_i)^{n_i}\phi^i_{j_i}(s_i)\right]\mc E_k\left(\frac{z^2}{\la}\right)\frac{z^{\sum_i(d_i-n_i)s_i}e^\la}{\la^{1+\sum_id_is_i}}\,{\rm d}\la\,{\rm d}s_1\dots{\rm d}s_h,
\]where $k=0,1,\dots,1-h+\sum_{i=1}^hn_i$, $\gm_i$ are parabolas encircling the poles of the factors $\Gm(s_i)^{n_i}$, the functions $\phi^i_{j_i}(s_i)$ are defined by
\begin{empheq}[left={\phi_{j_i}^{i}(s_i):=}\empheqlbrace]{align*}
\exp\left({2\pi\sqrt{-1}j_is_i}\right),\quad n_i\text{ even},\\
\exp\left({2\pi\sqrt{-1}j_is_i}+\pi\sqrt{-1}s_i\right),\quad n_i\text{ odd},
\end{empheq}
for any $h$-tuple $\bm j=(j_1,\dots, j_h)$ with $0\leq j_i\leq n_i-1$.
\vskip2mm

\noindent 1.5.\,\,The paper is organized as follows. 

In Section \ref{sec1}, we recall basic notions and results in Gromov--Witten theory, Frobenius manifolds and quantum cohomology.

In Section \ref{sec2}, we first introduce the quantum differential equation of a smooth projective variety and recall several results on the cyclic stratum and master functions from \cite{Cot22}. Subsequently, we overview several results in the literature on the quantum cohomology of projective bundles, focusing on quantum analogs of the classical Leray--Hirsch theorem. We also recall the definition of Borel $(\bm\al,\bm\bt)$-multitransform in the analytic setting. We introduce the sequence of $\mc E_k$-functions and describe their series expansions and integral representations. Then, we formulate the main results of the paper.

In Section \ref{sec3}, after introducing the notion of topological--enumerative solution of the {\it qDE} and of the closely related $J$-function, we recall the statement of the Elezi--Brown theorem \cite{Ele05,Ele07,Bro14}. Consequently, we prove the main result of the paper, Theorem \ref{mt2}.

Finally, in Section \ref{sec4} we exemplify our results on the specific case of the {\it qDE} of the blow-up of the projective plane at a point. We make explicit a base of solutions, obtained from a base of solutions for the {\it qDE} of $\mathbb P^1$ via the Borel multitransforms.
\vskip2mm
\noindent{\bf Acknowledgments. }
The author is thankful to U.\,Bruzzo, R. Conti, A.B.\,Cruzeiro, A.\,D'Agnolo, G. Degano, L.\,Fiorot, C.\,Florentino, D.\,Guzzetti, C.\,Hertling, M.\,Hien, P. Lorenzoni, D.\,Masoero, T.\,Monteiro Fernandes, G. Muratore, P.\,Polesello, L.\,Prelli, P.\,Rossi, A.T. Ricolfi, C. Sabbah, M.\,Smirnov, D.\,Van Straten, A.\,Varchenko, and J.C.\,Zambrini for several valuable discussions. This research was supported by the FCT projects PTDC/MAT-PUR/ 30234/2017, UIDB/00208/2020, 2021.01521.CEECIND, and 2022.03702.PTDC (GENIDE).

\section{Gromov--Witten invariants and quantum cohomology}\label{sec1}
\subsection{Notations}\label{secback}
Let $X$ be a smooth projective variety over $\C$ with vanishing odd cohomology, i.e.\,\,$H^{2k+1}(X,\C)=0$ for any $k\geq 0$. Fix a homogeneous basis $(T_0,\dots, T_n)$ of the $\C$-vector space $H^\bullet(X,\C)$, and denote by $\bm t=(t^0,\dots, t^n)$ the dual coordinates. Without loss of generality, we may assume $T_0=1\in H^0(X,\C)$. The cohomology $H^\bullet(X,\C)$ is equipped with a $\C$-algebra structure with respect to the (classical) cohomological product $\cup\colon H^\bullet(X,\C)\otimes H^\bullet(X,\C)\to H^\bullet(X,\C)$.
\vskip2mm
Consider the {\it Poincar\'e metric}\footnote{If $V$ is a finite dimensional $\C$-vector space, we denote by $\bigodot^kV$ the $k$-th symmetric tensor power of $V$.} $\eta\colon{\bigodot}^2 H^\bullet(X,\C)\to\C$, namely the non-degenerate symmetric $\C$-bilinear form on $H^\bullet(X,\C)$ defined by 
\[\eta(v_1,v_2):=\int_Xv_1\cup v_2,\quad v_1,v_2\in H^\bullet(X,\C).
\]For $\al,\bt=0,\dots, n$, we set $\eta_{\al\bt}:=\eta(T_\al,T_\bt)$, and we define $\eta:=\left(\eta_{\al\bt}\right)_{\al,\bt=0}^n$ to be the Gram matrix of $\eta$ with respect to the fixed basis $(T_0,\dots, T_n)$. The entries of the inverse Gram matrix $\eta^{-1}$ will be denoted by $\eta^{\al\bt}$, for $\al,\bt=0,\dots,n$. 
\vskip2mm
In all the paper, unless otherwise stated, the Einstein summation rule over repeated Greek indices will be used.
\subsection{Gromov--Witten invariants} Set $H_2(X,\Z)_{\rm tf}:=H_2(X,\Z)/{\rm torsion}$. 
\vskip2mm
Given $\bt\in H_2(X,\Z)_{\rm tf}$, denote by $\DM{0}{k}(X,\bt)$ the Kontsevich--Manin moduli stack of stable $k$-pointed rational maps of degree $\bt$ and with target $X$. This stack parametrizes isomorphism classes of triples $(C,\bm x, f)$ where 
\begin{enumerate}
\item $C$ is an algebraic curve of genus $0$ with at most nodal singularities,
\item $\bm x=(x_1,\dots, x_k)$ is a $k$-tuple of pairwise distinct points of the smooth locus of $C$,
\item $f\colon C\to X$ is a morphism such that $f_*[C]=\bt$,
\item a morphism $(C,\bm x, f)\to (C',\bm x', f')$ is the datum of a morphism $\phi\colon C\to C'$ such that $\phi(x_i)=x_i'$ for $i=1,\dots, k$ and $f=f'\circ \phi$,
\item the group of automorphisms of $(C,\bm x,f)$ is finite.
\end{enumerate}
The moduli space $\DM{0}{k}(X,\bt)$ is a proper Deligne--Mumford stack of virtual dimension
\[{\rm vir}\,\dim_\C \DM{0}{k}(X,\bt)=\dim_\C X-3+k+\int_\bt c_1(X).
\]

The  ({\it genus $0$}) {\it Gromov--Witten invariants} of $X$, and their {\it descendants}, are rational numbers $\langle\tau_{d_1}\gm_1,\dots, \tau_{d_k}\gm_k\rangle^X_{0,k,\bt}\in\Q$ defined via the intersection numbers of cycles on $\DM{0}{k}(X,\bt)$. More precisely, they are defined by the integrals
\beq\label{gw1}
\langle\tau_{d_1}\gm_1,\dots, \tau_{d_k}\gm_k\rangle^X_{0,k,\bt}:=\int_{\left[\DM{0}{k}(X,\bt)\right]^{\rm virt}}\prod_{i=1}^k{\rm ev}_i^*\gm_i\cup \psi_i^{d_i},
\eeq
where, for each $i=1,\dots, k$, we have
\begin{enumerate}
\item $\gm_i\in H^\bullet(X,\C)$,
\item $d_i\in \N$,
\item the $i$-th {\it evaluation morphism} ${\rm ev}_i\colon \DM{0}{k}(X,\bt)\to X$ acts by $(f,C,\bm x)\mapsto f(x_i)\in X$,
\item $\psi_i:=c_1(\mathscr L_i)$, where $\mathscr L_i\to \DM{0}{k}(X,\bt)$ is the $i$-th {\it tautological line bundle} whose fiber over $(C,\bm x, f)$ is $\mathscr L_i|_{(C,\bm x, f)}:=T^*_{x_i}C$.
\end{enumerate}
In the integral \eqref{gw1} the {\it virtual fundamental cycle} $\left[\DM{0}{k}(X,\bt)\right]^{\rm virt}$ is an element of the Chow ring ${\rm CH}_\bullet (\DM{0}{k}(X,\bt))$ with the expected dimension, that is
\[\left[\DM{0}{k}(X,\bt)\right]^{\rm virt}\in {\rm CH}_D (\DM{0}{k}(X,\bt)),\quad D={\rm vir}\,\dim_\C \DM{0}{k}(X,\bt).
\] See \cite{BF97} for its construction.

\subsection{Mori and K\"ahler cones} The stack $\DM{0}{k}(X,\bt)$ is non-empty only if $\bt$ is an element of the {\it Mori cone of effective curves}. This is the semigroup ${\rm NE}(X)$ of classes in the lattice $H_2(X,\Z)_{\rm tf}$ representable by algebraic curves, i.e.\,\,the classes $\bt=f_*[C]$ for some algebraic map $f\colon C\to X$, with $C$ a complete curve.

\begin{rem}\label{gwzero}
If the integrand \eqref{gw1} has a zero component in $H^{2D}(\DM{0}{k}(X,\bt))$, with $D={\rm vir}\,\dim_\C \DM{0}{k}(X,\bt)$, or if $\bt\not\in {\rm NE}(X)$, then $\langle\tau_{d_1}\gm_1,\dots, \tau_{d_k}\gm_k\rangle^X_{0,k,\bt}$ is set to 0.
\end{rem}

Inside $H^2(X,\R)\cap H^{1,1}(X,\C)$ we have a distinguished open convex cone, the {\it K\"ahler cone of $X$}, denoted by $\mc K_X$. This is defined as the set of all cohomology classes $[\om]\in H^{1,1}(X,\C)$ associated with any K\"ahler form $\om$ on $X$.
\vskip2mm
The Mori and K\"ahler cones are dual to each other, in the following sense.
\vskip2mm
By extension of scalars, we have an inclusion $i\colon H_2(X,\Z)_{\rm tf}\to H_2(X,\R)$. We will identify $H_2(X,\Z)_{\rm tf}$ with its image along $i$. 
\vskip2mm
Consider the smallest convex cone ${\rm NE}(X)_{\R}\subseteq H_2(X,\R)$ containing ${\rm NE}(X)$. 
\begin{prop}\label{mkcone}
The dual cone of the closure of ${\rm NE}(X)_\R$ equals the closure of $\mc K_X$.\qed
\end{prop}
Proposition \ref{mkcone} follows from 
\begin{itemize}
\item the identification of $\mc K_X$ with the cone of ample divisors of $X$, a consequence of Kodaira's characterization of positive classes (see e.g.\,\,\cite[Ch.\,5.3]{Huy04}\cite[Ch.\,7.2]{Voi02}),
\item and the Kleiman's ampleness criterion (see e.g.\,\,
\cite[Th.\,1.18]{KM98}\cite[Th.\,1.4.29]{Laz04}).
\end{itemize}
In particular, we have 
\beq\label{ample}
\int_\bt\om>0\qquad\text{for any}\quad\bt\in{\overline{\rm NE}(X)_\R}\setminus\{0\},\quad \om\in\mc K_X. 
\eeq
It is standard to denote by ${\rm Nef}(X)$ the dual cone of the closure $\overline{\rm NE}(X)_\R$, and to call it the {\it cone of nef divisors}. In this notation, the statement of Proposition \ref{mkcone} becomes ${\rm Nef}(X)=\overline{\mc K_X}$. We recommend the monograph \cite{Laz04} for a comprehensive account of the general theory of ample, nef, effective cones on smooth projective varieties.

\subsection{Novikov ring, and Gromov--Witten potentials} Fix a K\"ahler form $\om$ on $X$. We define the {\it Novikov ring} $\La_{X,\om}$ to be the ring of formal power series, in an indeterminate ${\bf Q}$, of the form
\[\sum_{\bt\in H_2(X,\Z)_{\rm tf}}a_\bt{\bf Q^\bt},\quad a_\bt\in\C,\quad{\rm card}\left\{\bt\colon a_\bt\neq 0\text{ and }\int_\bt\om<C\right\}<\infty,\quad \forall C\in\R.
\]
It is easy to see that the product of two such a series is well defined. In simple terms, $\La_{X,\om}$ can be intended as an ``upward'' completion of the group ring $\C[H_2(X,\Z)_{\rm tf}]$, by allowing sums with infinite terms in the direction of increasing values of the functional $\bt\mapsto \int_\bt\om$.
\vskip2mm
In what follows, it is convenient to associate with $X$ also a {\it big phase space} $\mc P_X$. This consists of an infinite product of countably many copies of $H^\bullet(X,\C)$, that is $\mc P_X:=H^\bullet(X,\C)^\N$. We will identify $H^\bullet(X,\C)$ with the $0$-th factor of $\mc P_X$, called the {\it small phase space}. Correspondingly to the choice of the base $(T_0,\dots, T_n)$, we denote by\footnote{The element $\tau_kT_\al$ is called the {\it descendant} of $T_\al$ with level $k$.} $(\tau_kT_0,\dots, \tau_kT_n)$ the basis of the $k$-th copy of $H^\bullet(X,\C)$ in $\mc P_X$. The coordinates of a point $\bm\gm\in \mc P_X$ with respect to $(\tau_kT_\al)_{\al,k}$ will be denoted by $\bm t^\bullet=\left(t^{\al,k}\right)_{\al,k}$. Hence, we will identify $t^\al\equiv t^{\al,0}$ for $\al=0,\dots, n$.
Instead of denoting by $\bm \gm=(t^{\al,k}\tau_kT_\al)_{\al,k}$ a generic point of $\mc P_X$ we will write this as a formal series
\[\bm\gm=\sum_{\al=0}^n\sum_{k=0}^\infty t^{\al,k}\tau_kT_\al.
\]
\vskip2mm
The ({\it genus 0}) {\it total descendant potential} of $X$ is the generating function $\mc F^X_0\in\La_{X,\om}[\![\bm t^\bullet]\!]$ of descendant Gromov--Witten invariants of $X$, defined by
\begin{multline*}
\mc F^X_0(\bm t^\bullet, {\bf Q}):=\sum_{k=0}^\infty\sum_{\bt\in H_2(X,\Z)_{\rm tf}}\frac{{\bf Q}^\bt}{k!}\langle\bm\gm,\dots,\bm\gm\rangle^X_{0,k,\bt}\\
=\sum_{k=0}^\infty\sum_{\bt}\sum_{\al_1,\dots,\al_k=0}^n\sum_{p_1,\dots,p_k=0}^\infty\frac{t^{\al_1,p_1}\dots t^{\al_k,p_k}}{k!}\langle\tau_{p_1}T_{\al_1}\dots\tau_{p_k}T_{\al_k}\rangle^X_{0,k,\bt}{\bf Q}^\bt.
\end{multline*}
The restriction of $\mc F^X_0$ to the small phase space (i.e.\,\,by setting $t^{\al,0}=t^{\al}$ and $t^{\al,p}=0$ for $p>0$) defines the ({\it genus 0}) {\it Gromov--Witten potential} of $X$, 
\[F^X_0(\bm t,{\bf Q}):=\sum_{k=0}^\infty\sum_\bt\sum_{\al_1,\dots, \al_k=0}^n\frac{t^{\al_1}\dots t^{\al_k}}{k!}\langle T_{\al_1},\dots, T_{\al_k}\rangle^X_{0,k,\bt}{\bf Q}^\bt,\qquad F^X_0\in\La_{X,\om}[\![\bm t]\!].
\]
\begin{rem}
The above formulas define two elements of $\La_{X,\om}[\![\bm t^\bullet]\!]$ (resp.\,$\La_{X,\om}[\![\bm t]\!]$) due to Remark \ref{gwzero} and ampleness condition \eqref{ample}. See also  \cite[Cor.\,1.19]{KM98}\cite[Ex.\,1.4.31]{Laz04}.
\end{rem}

\subsection{Quantum cohomology as $\La_{X,\om}$-formal Frobenius manifold} Consider
\begin{itemize}
\item the finite rank free $\La_{X,\om}$-module $H:=H^\bullet(X,\C)\otimes_\C\La_{X,\om}$,
\item the completed symmetric algebra $K:=\La_{X,\om}[\![H^T]\!]$ of the dual module $H^T$.
\end{itemize}
By extension of scalars, the $\C$-basis $(T_0,\dots, T_n)$ of $H^\bullet(X,\C)$ induces a basis of the module $H$. The algebra $K$ can be then identified with the formal power series ring $\La_{X,\om}[\![\bm t]\!]$. In what follows we set $H_K:=H\otimes_\C K$.
\vskip2mm
The formal spectrum $M:={\rm Spf}(K)$ is the formal scheme structure supported on a point (the ``origin'') with structure sheaf of formal functions $\Gm(M,\mathscr O_M)\cong K$, and with space of formal vector fields $\Gm(M,\mathscr T_M)\cong H_K$. Moreover, the Poincar\'e metric $\eta$ on $H^\bullet(X,\C)$ can be $\mathscr O_M$-bilinearly extended to a non-degenerate symmetric pairing $\eta\colon \mathscr T_M\otimes \mathscr T_M\to\mathscr O_M$.
\vskip2mm
The Gromov--Witten potential $F^X_0$ can be thus interpreted as a formal function on $M$. 
\begin{thm}\cite{KM94,Man99}
The function $F^X_0\in\mathscr O_M$ satisfies the following properties.
\begin{enumerate}
\item It is a quasi-homogeneous function. More precisely, if $E\in\mathscr T_M$ denotes the formal vector field on $M$
\beq\label{euler} 
E=c_1(X)+\sum_{\al=0}^n\left(1-\frac{1}{2}\deg T_\al\right)t^\al\frac{\der}{\der t^\al},
\eeq
then we have $EF^X_0=(3-\dim_\C X)F^X_0+Q$, where $Q$ is a quadratic polynomial in $\bm t$.
\item It is a solution of the Witten--Dijkgraaf--Verlinde--Verlinde (WDVV) equations, that is
\beq\label{wdvv}
\pushQED{\qed}
\sum_{\la,\nu}\frac{\der^3 F^X_0}{\der t^\al\der t^\bt\der t^\la}\eta^{\la\nu}\frac{\der^3 F^X_0}{\der t^\nu\der t^\gm\der t^\dl}=
\sum_{\la,\nu}\frac{\der^3 F^X_0}{\der t^\dl\der t^\bt\der t^\la}\eta^{\la\nu}\frac{\der^3 F^X_0}{\der t^\nu\der t^\gm\der t^\al}.\qedhere
\popQED
\eeq
\end{enumerate}
\end{thm}
The vector field $E$ of equation \eqref{euler} will be called {\it Euler vector field}.
\vskip2mm
Define the $K$-linear product $\circ$ of formal vector fields on $M$ via the formulas
\beq\label{Fr1}
T_\al\circ T_\bt:=\sum_\gm c^\gm_{\al\bt}T_\gm,\quad c^\gm_{\al\bt}:=\sum_\la \eta^{\gm\la}\frac{\der^3 F^X_0}{\der t^\al \der t^\bt \der t^\la},\quad \al,\bt,\gm=0,\dots, n.
\eeq
The product $\circ$ is associative, by WDVV equations \eqref{wdvv}, commutative, and compatible with the Poincar\'e pairing, as follows
\[\eta(v_1\circ v_2,v_3)=\eta(v_1,v_2\circ v_3),\quad \text{for any }v_1,v_2,v_3\in H_K.
\]Moreover, the vector field $e:=T_0$ is the unit for the $\circ$-product. This makes $(H_K,\circ,\eta,e)$ a Frobenius algebra.
\vskip2mm
The datum of $(M,F^X_0, \eta, e,E)$ defines a {\it formal Frobenius manifold} structure over the Novikov ring $\La_{X,\om}$, called {\it quantum cohomology} of $X$. See \cite{Man99}.

\subsection{Quantum cohomology as $\C$-formal Frobenius manifold}If we fix an arbitrary basis $(\bt_1,\dots, \bt_r)$ of $H_2(X,\Z)_{\rm tf}$, we can identify the Novikov ring $\La_{X,\om}$ with the ring of formal power series, in $r$ indeterminates ${ Q}_1,\dots, {Q}_r$ (where ${ Q}_i\equiv {\bf Q}^{\bt_i}$), of the form
\[\sum_{\bm d\in\Z^r}a_{\bm d}{\bf Q}^{\bm d},\quad {\bf Q}^{\bm d}:=\prod_{i=1}^r{Q}_i^{d_i},\quad {\rm card}\left\{\bm d\in\Z^r\colon a_{\bm d}\neq0,\int_{\sum_{i=1}^rd_i\bt_i}\om<C\right\}<\infty\quad \forall C\in\R.
\]Series as such are called {\it generalized Laurent series} in \cite{HS95}. 
\vskip1,5mm
Denote by $\mc B_X$ the set of points $\bm q\in (\C^*)^r$ such that for any $k\geq 3$, and any $\al_1,\dots,\al_k\in\{0,\dots,n\}$, the generalized Laurent series 
\beq\label{sumA}
\sum_{\bm d\in\Z^r}\langle T_{\al_1},\dots, T_{\al_k}\rangle^X_{0,k,\sum_{i=1}^rd_i\bt_i}{\bf Q}^{\bm d}
\eeq
is convergent at ${\bf Q}=\bm q$.
\begin{rem}
By multi-linearity of Gromov--Witten invariants, the sum \eqref{sumA} is convergent if and only if all the series $$\sum_{\bm d\in\Z^r}\langle v_{\al_1},\dots, v_{\al_k}\rangle^X_{0,k,\sum_{i=1}^rd_i\bt_i}{\bf Q}^{\bm d}$$ are convergent at ${\bf Q}=\bm q$ for any $k\geq 3$, any $\al_1,\dots,\al_k\in\{0,\dots,n\}$, and any choice of vectors $v_{\al_1},\dots, v_{\al_k}\in H^\bullet(X,\C)$.
\end{rem}
\vskip1,5mm
\noindent{\bf Assumption A:} The set $\mc B_X$ is non-empty.
\vskip1,5mm
\begin{prop}\label{assA}
Assumption A holds if and only if $\mc B_X=(\C^*)^r$.
\end{prop}
\proof
Let $\bt_1^\vee,\dots, \bt_r^\vee\in H^2(X,\C)$ be the dual classes of $\bt_1,\dots, \bt_r$, so that $\int_{\bt_i}\bt_j^\vee=\dl_{ij}$. Fix $\bm q_0,\bm q_1\in (\C^*)^r$.
By the divisor axiom of Gromov--Witten invariants, we have\footnote{We use the multi-index notation ${\bm q}^{\bm d}=(q_1^{d_1},\dots, q_r^{d_r})$, for $r$-tuples $\bm q,\bm d$.}
\begin{multline*}\sum_{\bm d\in\Z^r}\langle T_{\al_1},\dots, T_{\al_k}\rangle^X_{0,k,\sum_{i=1}^rd_i\bt_i}{\bm q_1}^{\bm d}=\sum_{\bm d\in\Z^r}\langle T_{\al_1},\dots, T_{\al_k}\rangle^X_{0,k,\sum_{i=1}^rd_i\bt_i}\left(\bm q_0^{-1}\bm q_1\right)^{\bm d}{\bm q_0}^{\bm d}\\=\sum_{\bm d\in\Z^r}\langle T_{\al_1},\dots, T_{\al_k}\rangle^X_{0,k,\sum_{i=1}^rd_i\bt_i}{\prod_{i=1}^r\exp \left(d_i\log \frac{q_{1,i}}{q_{0,i}}\right)}\bm q_0^{\bm d}\\
=\sum_{\bm d\in\Z^r}\langle \log \frac{q_{1,1}}{q_{0,1}}\,\bt_1^\vee, \dots, \log \frac{q_{1,r}}{q_{0,r}}\bt_r^\vee,T_{\al_1},\dots, T_{\al_k}\rangle^X_{0,r+k,\sum_{i=1}^rd_i\bt_i}\bm q_0^{\bm d}.
\end{multline*}
Hence, $\bm q_1\in \mc B_X$ if and only if $\bm q_0\in \mc B_X$.
\endproof
\begin{prop}
If $X$ is Fano, each series \eqref{sumA} is finite. Hence, Assumption A holds true.
\end{prop}
\proof
The invariant $\langle T_{\al_1},\dots, T_{\al_k}\rangle^X_{0,k,\sum_{i=1}^rd_i\bt_i}$ is nonzero only if $\sum_{i=1}^rd_i\bt_i\in{\rm NE}(X)$, and we have the dimensional constraint $\sum_{i=1}^k\deg T_{\al_i}=\int_{\sum_{i=1}^rd_i\bt_i} c_1(X)+\dim_\C X+k-3$. If $X$ is Fano, then $c_1(X)$ is ample, so that $\int_\bt c_1(X)>0$ if $\bt \in {\rm NE}(X)$, by \eqref{ample}. As a consequence of Kleiman ampleness criterion \cite[Cor.\,1.19]{KM98}, there exist only finitely many tuples $\bm d\in\Z^r$ for which the dimensional constraint holds.
\endproof

Whenever Assumption A holds, the specialization $F^X_0|_{{\bf Q}=\bm q}$, with $\bm q\in \mc B_X$, is a well-defined formal power series with complex coefficients. A formal Frobenius manifold over $\C$,
\[\left({\rm Spf}\,\C[\![H^\bullet(X,\C)^T]\!],\quad F^X_0|_{{\bf Q}=\bm q},\quad \eta,\quad e,\quad E\right),
\]is defined for each $\bm q\in \mc B_X$, by specializing the Gromov--Witten potential $F^X_0$ at ${\bf Q}=\bm q$. We call such a formal Frobenius manifold {\it quantum cohomology of $X$ specialized at ${\bf Q}=\bm q$}.

\subsection{Quantum cohomology as Dubrovin--Frobenius manifold} In this paper we will consider the case of {\it convergent} Gromov--Witten potentials only. This will allow us to promote the formal Frobenius manifold structures to analytic ones. 
\vskip2mm
Consider a variety $X$ for which Assumption A holds, so that $F^X_0|_{{\bf Q}=\bm q}\in\C[\![\bm t]\!]$ for any $\bm q\in\mc B_X$. Define $\mc B_X'\subseteq \mc B_X$ to be the set of points $\bm q\in\mc B_X$ such that $F^X_0|_{{\bf Q}=\bm q}\in\C\{\bm t\}$.
\vskip2mm
\noindent {\bf Assumption B: }The set $\mc B_X'$ is non-empty.
\begin{prop}
Assumption B holds if and only if $\mc B'_X=\mc B_X=(\C^*)^r$.
\end{prop}
\proof
Let $\bm q_0,\bm q_1\in\mc B_X$, and set $F_0,F_1\in\C[\![\bm t]\!]$ be the formal power series defined by the specializations $F_{(0)}:=F^X_0|_{{\bf Q}=\bm q_0}$ and $F_{(1)}:=F^X_0|_{{\bf Q}=\bm q_1}$. The same computations in the proof of Proposition \ref{assA}, invoking the divisor axiom, show that $F_{(1)}$ is a derivative of $F_{(0)}$. More precisely, if $V_i$ denotes the vector field corresponding to the cohomology class $\log \frac{q_{1,i}}{q_{0,i}}\,\bt_i^\vee$ for $i=1,\dots, r$, we have $F_{(1)}=V_1\dots V_rF_{(0)}$. Hence $F_{(0)}\in \C\{\bm t\}$ if and only if $F_{(1)}\in \C\{\bm t\}$.
\endproof

The following result can be useful to deduce the validity of Assumption B.

\begin{thm}\cite{Cot21a,Cot21b}
Let Assumption A hold. Assume that the quantum cohomology of $X$ specialized at some point ${\bf Q}=\bm q_o$, with $\bm q_o\in\mc B_X$, is formally semisimple, that is the $\C$-algebra $(H^\bullet(X,\C),\circ_{\bm q_o,0})$ with structure constants $\sum_\dl\eta^{\la\dl}\der^3_{t^\al t^\bt t^\dl}F^X_0(\bm t,\bm q_o)|_{\bm t=0}$ is isomorphic to $\C^{n+1}$. Then Assumption B holds.\qed  
\end{thm}
Under the validity Assumption B, the series $F^X_0|_{{\bf Q}=\bm q}$ have the same domain of convergence $M\subseteq \C^{n+1}$ for any $\bm q\in (\C^*)^r$. Without loss of generality, in what follows we will consider the specialization of $F^X_0$ at ${\bf Q}=\bm 1=(1,1,\dots,1)$.
\vskip2mm
Let $TM$ (resp.\,\,$T^*M$) the holomorphic tangent (resp.\,\,cotangent) bundles of $M$. 
At each $p\in M$, we have a canonical identification of vector spaces $T_pM\cong H^\bullet(X,\C)$, via the map $\frac{\der}{\der t^\al}\mapsto T_\al$, for $\al=0,\dots, n$. 
Via this identification, the Poincar\'e metric defines a holomorphic section $\eta\in \Gm\left(\bigodot^2T^*M\right)$. This is a symmetric non-degenerate $\mc O_M$-bilinear 2-form, for simplicity called {\it metric}, with {\it flat} Levi--Civita connection $\nabla$.

Consider the holomorphic $(1,2)$-tensor $c\in \Gm(TM\otimes \bigodot^2 T^*M)$ with components $$c^\al_{\bt\gm}:=\sum_\la \eta^{\al\la}\nabla_\la\nabla_\bt\nabla_\gm F^X_0|_{{\bf Q}=\bm 1}=\sum_\la\eta^{\al\la}\frac{\der^3}{\der t^\la\der t^\bt \der t^\gm}F^X_0|_{{\bf Q}=\bm 1},$$
together with the two holomorphic vector fields $e,E\in \Gm(TM)$ defined by the constant field $e=T_0\in H^\bullet(X,\C)$ and equation \eqref{euler}, respectively.

The product defined in equation \eqref{Fr1} is now convergent: each tangent space $T_pM$ is equipped with a Frobenius algebra structure, holomorphically depending on the point $p$. The tangent vector $e_p$ equals the unit of the Frobenius algebra at $p\in M$. For this reason $e$ is called {\it unit vector field}. 
The {\it Euler vector field} $E$ satisfies the following identities:
\[\frak L_Ec=c,\quad \frak L_E\eta=(2-\dim_\C X)\eta.
\]The datum of $(M,\eta,c,e,E)$ defines a Dubrovin--Frobenius manifold structure \cite{Dub96,Dub99}.

\section{QDEs, Borel multitransforms, $\mc E_k$ functions, and main results }\label{sec2}

\subsection{Extended deformed connection}
Consider the Dubrovin--Frobenius manifold $M$ defined in the previous section. Let us introduce two distinguished holomorphic $(1,1)$-tensors on $M$.
\begin{defn}
The \emph{grading operator} $\bm \mu\in \Gm({\rm End}\,TM)$ is the tensor  defined by
\beq
\bm\mu(v):=\frac{2-\dim_\mathbb CX}{2}v-\nabla_vE,\quad v\in\Gamma(TM).
\eeq
The operator $\bm {\mc U}\in \Gm({\rm End}\,TM)$ is the tensor  defined by
\beq
\bm{\mc U}(v):=E\circ v,\quad v\in\Gamma(TM).
\eeq
We will denote by $\mu$ and $\mc U$, respectively, the matrices of components of the tensors $\bm\mu$ and $\bm{\mc U}$ in the $\nabla$-flat coordinates $\bm t$.
\end{defn}
\begin{lem}\label{lemsymUmu}
For any vector fields $v_1,v_2\in\Gm(TM)$, we have
\[\eta\left({\bm{\mc U}}(v_1),v_2\right)=\eta\left(v_1,{\bm {\mc U}}(v_2)\right),\qquad \eta\left({\bm{\mu}}(v_1),v_2\right)=-\eta\left(v_1,{\bm {\mu}}(v_2)\right).
\]
\end{lem}
\proof
The result directly follows from the definitions.
\endproof
\vskip2mm
Consider a punctured complex line $\C^*$, with global coordinate $z$, and introduce the extended space $\widehat M:=\C^*\times M$. Given the canonical projection $\pi\colon\widehat M\to M$, we can consider the pull-back bundle $\pi^*TM$. All the tensors $\eta,c, E, \bm\mu,\bm{\mathcal U}$ can be lifted to $\pi^*TM$, and their lifts will be denoted by the same symbols. Moreover, the Levi--Civita connection $\nabla$ can be uniquely lifted on $\pi^*TM$ in such a way that
\[\nabla_\frac{\partial}{\partial z}v=0\quad \text{for any }v\in\pi^{-1}\mathscr T_M,
\]where $\pi^{-1}\mathscr T_M$ the sheaf of sections of $\pi^*TM$ constant on the fibers of $\pi$.

\begin{defn}
The \emph{extended deformed connection} is the connection $\widehat \nabla$ on the bundle $\pi^*TM$ defined by
\begin{align}
\widehat\nabla_{w}v&=\nabla_wv+z\cdot w*v,\\
\widehat\nabla_{\frac{\partial}{\partial z}}v&=\nabla_{\partial_z}v+\mathcal U(v)-\frac{1}{z}\mu(v),
\end{align}
for $v,w\in \Gamma(\pi^*T\Omega)$.
\end{defn}

\begin{thm}[\cite{Dub96,Dub99}]
The connection $\widehat\nabla$ is flat.\qed
\end{thm}

\subsection{The quantum differential equation} The connection $\widehat\nabla$ induces a flat connection on $\pi^*T^*M$. Let $\varpi\in\Gamma(\pi^*T^*M)$ be a flat section. Consider the corresponding vector field $\sig\in\Gamma(\pi^*TM)$ via musical isomorphism\footnote{The datum of a non-degenerate symmetric bilinear form $\eta$ defines an isomorphism between tangent and cotangent bundles of $M$. A vector field $v\in\Gm(TM)$ and a 1-form $\xi\in\Gm(T^*M)$ are related to each other via these isomorphisms if and only if $\xi(w)=\eta(v,w)$ for any $w\in\Gm(TM)$. In components, these isomorphisms act as the operations of ``lowering'' and ``raising'' the indices. For this reason, these isomorphisms are oftentimes called {\it musical}, in analogy with the $\flat$ (flat) and $\sharp$ (sharp) operations on the musical notes on a staff: these respectively make a musical note ``lower'' or ``higher'' in pitch by a semitone. }, i.e. such that $\varpi(v)=\eta(\sig,v)$ for all $v\in\Gamma(\pi^*TM)$.

The vector field $\sig$ satisfies the following system\footnote{We consider the joint system \eqref{eq1}, \eqref{qde} in matrix notations ($\sig$ is a column vector whose entries are the components $\sig^\alpha(\bm t,z)$ wrt $\frac{\partial}{\partial t^\alpha}$). Bases of solutions are arranged in invertible $n\times n$-matrices, called \emph{fundamental systems of solutions}.} of equations
 \begin{align}
 \label{eq1}
 \frac{\partial}{\partial t^\alpha}\sig&=z\,\mathcal C_\alpha\sig,\quad \alpha=0,\dots, n,\\
\label{qde}
  \frac{\partial}{\partial z}\sig&=\left(\mathcal U+\frac{1}{z}\mu\right)\sig.
 \end{align}
Here $\mathcal C_\alpha$ is the $(1,1)$-tensor defined by $(\mathcal C_\alpha)^\beta_\gamma:=c^\beta_{\alpha\gamma}$. 
\begin{defn}
The \emph{quantum differential equation} ($qDE$) of $X$ is the differential equation \eqref{qde}. 
\end{defn}

Equivalently, by invoking the (skew)symmetry properties of Lemma \ref{lemsymUmu}, the differential system \eqref{eq1}, \eqref{qde} can be written in terms of the unkown covector field $\varpi$: 
\begin{align}
\label{eq1.2}
\frac{\der}{\der t^\al}\varpi&=z\,\mc C_\al^T\varpi,\quad \al=0,\dots,n,\\
\label{qde.2}
\frac{\der}{\der z}\varpi&=\left(\mc U^T-\frac{1}{z}\mu\right)\varpi.
\end{align}
Here $\varpi:=\eta\sig$ is a column vector with components $(\varpi_\al)_{\al=0}^n$.

\subsection{Cyclic stratum, $\mc A_\La$-stratum, master functions}\label{scyc} For any $k\in\N$, introduce a vector field $e_k\in\Gm(\pi^*TM)$, defined by the iterated covariant derivative
\beq
e_k:=\edc_{\frac{\der}{\der z}}^ke.
\eeq
\begin{defn}
The {\it cyclic stratum} $\widehat M^{\rm cyc}$ is defined to be the maximal open subset $U$ of $\widehat M$ such that the bundle $\pi^*TM|_U$ is trivial and the collection of sections $(e_k|_{U})_{k=0}^{n}$ defines a basis of each fiber. On ${\widehat M}^{\rm cyc}$ we will also introduce the dual coframe $(\omega_j)_{j=0}^{n}$, by imposing
\beq
\label{eom}
\langle\omega_j,e_k\rangle=\delta_{jk}.
\eeq
The frame $(e_k)_{k=0}^{n}$ will be called \emph{cyclic frame}, and its dual $(\om_j)_{j=0}^{n}$ \emph{cyclic coframe}.
\end{defn}

\begin{defn}The $\La$-{\it matrix} is the matrix-valued function $\La=(\La_{i\al}(z,p))$, holomorphic on ${\widehat M}^{\rm cyc}$, defined by the equation
\beq\label{Lam}
\frac{\der}{\der t^\al}=\sum_{i=0}^{n}\La_{i\al}e_i,\quad \al=0,\dots,n.
\eeq
\end{defn}

\begin{thm}\label{strLa}\cite[Th.\,2.20]{Cot22}
The function $\det\La$ is a meromorphic function on $\Pb^1\times M$  of the form
\[\det\La(z,p)=\frac{z^{\binom{n}{2}}}{z^{\binom{n}{2}}A_0(p)+\dots+A_{\binom{n}{2}}(p)},
\]
where $A_0,\dots, A_{\binom{n}{2}}$ are holomorphic functions on $M$. Moreover, if $n>1$ and if the eigenvalues of the grading operator $\bm\mu$ are not pairwise distinct, then the function $A_{\binom{n}{2}}$ is identically zero. \qed
\end{thm}

If we set $$\bar n:=\min\{j\in\N\colon A_{h}(p)=0\quad\forall p\in M,\,\forall h>j\},$$ the function $\det\La$ takes the form
\[\det \La=\frac{z^{\bar n}}{z^{\bar n} A_0(p)+z^{\bar n-1}A_1(p)\dots+A_{\bar n}(p)}.
\]
\begin{defn}We define the set $\mc A_{\La}\subseteq M$ to be the set 
\[\mc A_\La:=\left\{p\in M\colon\quad A_0(p)=\dots=A_{\bar n}(p)=0\right\}.
\]
\end{defn}
\vskip2mm
Consider now the system \eqref{eq1.2}, \eqref{qde.2}. We have the following result.

\begin{thm} \cite[Th.\,2.29]{Cot22}\label{thsqde}
The matrix differential equation \eqref{qde.2}, specialized at a point $p\in M\setminus\mc A_{\La}$, can be reduced to a single scalar differential equation of order $n+1$ in the unknown function $\varpi_0$. The scalar differential equation admits at most $\binom{n}{2}$ apparent singularities. \qed
\end{thm}
The scalar differential equation above will be called the {\it master differential equation} of $X$.
\vskip2mm
We can more explicitly describe the master differential equation.
At points $(z,p)\in \widehat M^{\rm cyc}$, introduce the column vector $\overline{\varpi}$ defined by
\[\overline{\varpi}:=\left(\La^{-1}\right)^T\varpi,
\]where $\La$ is the matrix defined as in \eqref{Lam}. The entries of $\overline{\varpi}$ are the components $\overline{\varpi}_j$ of the $\edc$-flat covector $\varpi$ with respect to the cyclic coframe $(\om_j)_{j=0}^n$. The vector $\overline\varpi$ satisfies the differential system

\begin{align}
\label{eq1.3}
\frac{\der \overline\varpi}{\der t^\al}&=\left(z\left(\La^{-1}\right)^T\mc C_\al^T\La^T+\frac{\der\left(\La^{-1}\right)^T}{\der t^\al}\La^T\right)\overline\varpi,\\
\label{qde.3}
\frac{\der \overline\varpi}{\der z}&=\left(\left(\La^{-1}\right)^T\mc U^T\La^T-\frac{1}{z}\left(\La^{-1}\right)^T\mu\La^T+\frac{\der\left(\La^{-1}\right)^T}{\der z}\La^T\right)\overline\varpi.
\end{align}

\begin{thm}\cite[Cor.\,2.27]{Cot22}
The system of differential equations \eqref{qde.3} is the companion system of the master differential equation of $X$. \qed
\end{thm}
\begin{rem}
Since $e_0=e$, we have $\varpi_0=\overline{\varpi}_0$.
\end{rem}

\begin{defn}\label{defmastfun} Fix $p\in M$, consider the system of differential equations \eqref{qde.2} specialized at $p$, and let $\mc X_p$ be the $\C$-vector space of its solutions. Let $\nu_p\colon\mc X_p\to\mc O\left(\widetilde{\C^*}\right)$ be the morphism defined by
\[\varpi\mapsto\Phi_\varpi(z),\quad \Phi_\varpi(z):=z^{-\frac{\dim X}{2}}\langle \varpi(z,p),e(p)\rangle.
\]Set $\mc S_p(M):={\rm im\,}\nu_p$. Elements of $\mc S_p(M)$ will be called {\it master functions} of $X$ at $p$.
\end{defn}

By Theorem \ref{thsqde}, the morphism $\nu_p$ is injective at points $p\in M\setminus\mc A_\La$, see \cite[Th.\,2.31]{Cot22}.

\subsection{Projective bundles: classical and quantum aspects}\label{secpb}
Let $X$ be a smooth projective variety over $\C$, and $V\to X$ a holomorphic vector bundle on $X$ of rank $r+1$. For short denote by $P$ the total space of the projective bundle\footnote{The fiber $P_x$ over $x\in X$ parametrizes one-dimensional subspaces of $V_x$. The reader should be aware that a conflicting definition is also in use. Some authors, following Grothendieck, define the projectivization $\Pb (V)$ by parametrizing hyperplanes (or 1-quotients) of fibers of $V$.} $\pi\colon \Pb (V)\to X$, and set 
\[\xi:=c_1\left(\mc O_P(1)\right),\quad \text{where }\mc O_P(-1)\text{ is the tautological line bundle on }P.
\] 
The classical cohomology $H^\bullet(P,\C)$ is an algebra over $H^\bullet(X,\C)$ via pullback. More precisely, by the classical Leray--Hirsch Theorem, the pull-back map $\pi^*\colon H^\bullet(X,\C)\to H^\bullet(P,\C)$ is a monomorphism of rings, and via this map the algebra $H^\bullet(P,\C)$ admits the following presentation
\[H^\bullet(P,\C)\cong H^\bullet(X,\C)[\xi]\Big/\!\!\left(\xi^{r+1}+c_1(V)\xi^{r}+\dots+c_{r+1}(V)\right).
\]Moreover, the $\C$-linear map $H^\bullet(X,\C)^{\oplus (r+1)}\to H^\bullet(P,\C)$ defined by $(\al_0,\dots,\al_r)\mapsto \sum_i\xi^i\pi^*\al_i$ is an isomorphism of $\C$-vector spaces, so that
\[H^\bullet(P,\C)\cong \bigoplus_{i=0}^r\xi^i H^\bullet(X,\C)\quad\text{ as $\C$-vector spaces.}
\]See e.g.\,\,\cite[Sec.\,\,11.2]{EH16}.
\begin{rem}\label{rklines}
Two projective bundles $\Pb(V)$ and $\Pb(V')$ are isomorphic as $X$-schemes if and only if $V'\cong V\otimes L$ for some line bundle $L\to X$. In such a case, we have $\mc O_{P'}(-1)\cong \pi^*L\otimes \mc O_{P}(-1)$, so that $\xi'=\xi-c_1(L)$. See \cite[Prop.\,\,11.3]{EH16}. Consequently, in the case of a split bundle $V=\bigoplus_{j=0}^rL_j$, without loss of generality we can assume that $L_0=\mc O_X$.
\end{rem}
We say that $V\to X$ is a {\it Fano bundle} if the projectivization $P=\Pb(V)$ is a Fano manifold. The existence of a Fano bundle of $X$ automatically implies that $X$ itself is Fano \cite[Th.\,1.6]{SW90-1}. Fano bundles have been extensively studied, and even completely classified in several cases:
\begin{enumerate}
\item In \cite{SW90-1,SW90-2,SW90-3,SSW91} a complete classification of rank two Fano bundles up to dimension 3 is given.
\item In \cite{APW94} it is shown that a rank two Fano bundle on $\Pb^n$, with $n\geq 4$, and a quadric $Q_n$, with $n\geq 6$, splits into a direct sum of line bundles. On $Q_4$ and $Q_5$ there are some exceptions: two spinor bundles and a 7-dimensional family of stable bundles on $Q_4$, and Cayley bundles on $Q_5$ \cite{Ott90}.
\item In \cite{MOS12a,MOS12b}, rank two Fano bundles on projective Grassmannians $\mathbb G(1,n)$, with $n\geq 4$, have been classified. It is shown that the only non-split rank two Fano bundles are, up to twists, the universal quotient bundles $\mc Q\to\mathbb G(1,n)$. Subsequently, in \cite{MOS14}, a classification of rank two Fano bundles on manifolds $X$ with $H^2(X,\Z)\cong H^4(X,\Z)\cong \Z$ is given.
\end{enumerate}
If $E\to X$ is a complex vector bundle, let $c(E)$ be the total Chern class $c(E)=\sum_{j\geq 0}c_j(E)$. The next result shows how to compute the whole Chern class $c(P)$ 
in terms of $c(X)$, $c(V)$, and $\xi$.
\begin{lem}\label{lccp}
We have $c(P)=c(\pi^*V\otimes \mc O_P(1))\pi^*c(X)$. In particular, we have
\[c_1(P)=\pi^*c_1(X)+\pi^*c_1(V)+(r+1)\xi.
\]
\end{lem}
\proof The first claim follows from the relative Euler exact sequence 
\[0\to \mc O_{P}\to\pi^*V\otimes \mc O_P(1)\to T_{P/X}\to 0,
\]and the exact sequence
\[
0\to T_{P/X}\to TP\to \pi^*TX\to 0.
\]The $k$-th Chern class of $\pi^*V\otimes \mc O_P(1)$ can be computed by the formula
\[
\pushQED{\qed}
c_k(\pi^*V\otimes \mc O_P(1))=\sum_{j=0}^{k}\binom{r+1-j}{k-j}c_j(\pi^*V)\xi^{k-j}, \quad k=0,\dots,r+1.\qedhere
\popQED
\]

In \cite{MP06}, D.\,Maulik and R.\,Pandharipande developed a quantum Leray--Hirsch theorem for $\Pb^1$-bundles $P=\Pb(\mc O_X\oplus L)$ over a smooth base $X$. In particular, they showed that the Gromov--Witten theory of $P$ can be uniquely and effectively reconstructed from the Gromov--Witten theory of $X$ and the Chern class $c_1(L)$.

Contemporarily, in \cite{Ele05} A.\,Elezi formulated a conjectural quantum Leray--Hirsch result, which allows the reconstruction of the small $J$-function of $P$ starting from the $J$-function of $X$, in the case $P$ is the projectivization of a split Fano bundle on $X$. The validity of this conjecture was first proved in several cases in \cite{Ele07}. Subsequently, in \cite{Bro14} J.\,Brown proved Elezi conjecture as a special case of a more general result, allowing the reconstruction of the $J$-function of the total space of a toric fibration $F\to X$ (constructed from a split vector bundle $E\to X$) from the $J$-function of the base $X$. 

The result of Brown was the starting point for a big project of Y.-P.\,Lee, H.-W.\,Lin and C.-L.\,Wang, developed in the series of papers \cite{LLW10,LLW16a,LLW16b}. The original Brown's theorem was formulated in terms of Givental's Lagrangian cones formalism \cite{Giv04,CG07}. Lee, Lin and Wang first showed how to explicitly reconstruct the $J$-function of a toric fibration in terms of {\it generalized mirror transformations} and {\it Birkhoff factorizations}. Then, they used these results to show the invariance of quantum cohomology under ordinary flops. Moreover, in a subsequent paper \cite{LLQW16} joint with W.\,Qu, they showed how it is possible to remove any splitting assumption in the quantum Leray--Hirsch theorem ({\it quantum splitting principle}).

In \cite{Fan21}, H.\,Fan recently proved that the Gromov--Witten theory of a projective bundle $\Pb(V)\to X$ is uniquely determined by the Gromov--Witten theory of $X$ and the Chern class $c(V)$.

In the very recent preprint \cite{Kot22}, it is shown that if Assumption $B$ holds for the quantum cohomology of $X$, then the same holds for the quantum cohomology of $P$. In particular, its quantum cohomology can be equipped with a Dubrovin--Frobenius manifold structure. 

\subsection{Analytic Borel multitransform} Let $(A,+,\cdot,1_A)$ be an associative, commutative, unital finite dimensional $\C$-algebra.
\vskip2mm
Let $h\in\N^*$, and $\bm\al,\bm\bt\in(\C^*)^h$. Given a $h$-tuple $(\Phi_1,\dots,\Phi_h)$ of $A$-valued holomorphic functions defined on $\C\setminus ]-\infty;0]$, we define their {\it Borel }$(\bm\al,\bm\bt)$-{\it multitransform} by the integral (if convergent)
\[\mathscr B_{\bm\al,\bm\bt}\left[\Phi_1,\dots,\Phi_h\right](z):=\frac{1}{2\pi\sqrt{-1}}\int_{\frak H}\prod_{j=1}^h\Phi_j\left(z^\frac{1}{\al_j\bt_j}\la^{-\bt_j}\right)e^\la\frac{{\rm d}\la}{\la},
\]where $\frak H$ is a Hankel-type contour of integration, see Figure \ref{gammahankel}.

\begin{prop}{\cite[Prop.\,6.12]{Cot22}}
Let $(e_1,\dots, e_n)$ be a $\C$-basis of $A$. Given a tuple $(\Phi_1,\dots,\Phi_h)$ of $A$-valued functions, introduce the corresponding $\C$-valued component functions $\Phi_j^i$, with $i=1,\dots, n$, $j=1,\dots, h$, such that\,\,$\Phi_j=\sum_{i=1}^n\Phi_j^ie_i$. The components of $\mathscr B_{\bm\al,\bm\bt}\left[\Phi_1,\dots,\Phi_h\right]$ are $\C$-linear combinations of the $h\cdot n$ $\C$-valued functions $\mathscr B_{\bm\al,\bm\bt}\left[\Phi_1^{i_1},\dots,\Phi_h^{i_h}\right]$, where $(i_1,\dots,i_h)\in\{1,\dots,n\}^{\times h}$.\qed
\end{prop}
\subsection{The functions $\mc E_k$}\label{secek}

The {\it B\"ohmer--Tricomi normalized Gamma function} $\gm^*$ is the entire function on $\C^2$, defined by the integral
\[\gamma^*(s,z):=\frac{1}{\Gm(s)}\int_0^1t^{s-1}e^{-zt}{\rm d} t,\qquad {\rm Re}\,s>0,\quad z\in\C.
\]This function was originally introduced by P.E.\,B\"ohmer in \cite[pp.\,124--125]{B\"oh39}, and subsequently studied in \cite{Tri50}. It is closely related to the {\it upper} and {\it lower incomplete Gamma functions} $\Gm$ and $\gm$, introduced in 1877 by F.E.\,Prym \cite{Pry77}. These are the holomorphic functions defined on $\C\times \widetilde{\C^*}$ whose general values are respectively defined by 
\[\Gm(s,z):=\int_z^\infty t^{s-1}e^{-t}{\rm d}t,\qquad \gm(s,z):=\int_0^z t^{s-1}e^{-t}{\rm d}t,\qquad{\rm Re}\,s>0,\quad z\in\C,
\]without restrictions on the integration paths. We have
\[\gm^*(s,z)=z^{-s}\frac{\gm(s,z)}{\Gm(s)}=z^{-s}\left(1-\frac{\Gm(s,z)}{\Gm(s)}\right).
\]For general properties of the incomplete Gamma functions, see e.g.\,\,\cite[Kap.\,II, XV, XXI]{Nie06} \cite[Kap.\,V]{B\"oh39} \cite[Ch.\,IX]{Erd53} \cite[\S 6.5]{AS64} \cite[Ch.\,8]{OLBC10}. 
\vskip2mm
Let $\mc E$ be the holomorphic function on $\C\times \widetilde {\C^*}$ defined by 
\beq\mc E(s,z):=z^s e^z\gm^*(s,z)=e^z\left(1-\frac{\Gm(s,z)}{\Gm(s)}\right).
\eeq

\begin{prop}\label{expe}
For any $(s,z)\in\C\times \widetilde {\C^*}$, we have
\[\mc E(s,z)=\sum_{k=0}^\infty\frac{z^{k+s}}{\Gm(1+k+s)}.
\]
\end{prop}
\proof
We have $\gm^*(s,z)=e^{-z}\sum_{k=0}^\infty\frac{z^{k}}{\Gm(1+k+s)}$, see \cite[Ch.\,8]{OLBC10}. 
The result follows.
\endproof

\begin{defn}
For any $k\in\N$, we define the holomorphic function $\mc E_k\in\mc O(\widetilde {\C^*})$ by the iterated partial derivative
\beq
\mc E_k(z):=\left.\frac{\der^k}{\der s^k}\mc E(s,z)\right|_{s=0},\qquad z\in\widetilde{\C^*}.
\eeq
\end{defn}
Introduce the following notations:
\begin{itemize}
\item $\dl_{ij}$ denotes the Kronecker symbol,
\item $\cem:=\lim_n\left(\sum_{k=1}^n\frac{1}{k}-\log n\right)$ denotes the Euler--Mascheroni constant,
\item $\zeta(s)$ denotes the Riemann zeta function,
\item $B_n(x_1,\dots,x_n)$ denotes the $n$-th complete Bell polynomial, recursively defined by
\beq\label{bellp}
B_0:=1,\quad B_{n+1}(x_1,\dots, x_{n+1}):=\sum_{k=0}^n\binom{n}{k}B_{n-k}(x_1,\dots, x_{n-k})x_{k+1},\quad n\geq 0,
\eeq
\item $(\bt_n)_{n\in\N}$ is the sequence of real numbers defined by
\[\bt_n:=B_n(0!\cem,\,-1!\zeta(2),\,2!\zeta(3),\,-3!\zeta(4),\,\dots,\, (-1)^{n+1}(n-1)!\zeta(n)),
\]
\item $G_{p,q}^{m,n}\left(z\left|\begin{array}{c}
a_1,\dots, a_p\\
b_1,\dots,b_q
\end{array}\right.\right)$ denotes the Meijer $G$-function defined by the integral
\[
\frac{1}{2\pi\sqrt{-1}}\int_L\frac{\prod_{i=1}^n\Gm(1-a_i+s)\,\prod_{i=1}^m\Gm(b_i-s)}{\prod_{i=n+1}^p\Gm(a_i-s)\,\prod_{i=m+1}^q\Gm(1-b_i+s)}z^{s}{\rm d}s,
\]for a suitable integration path $L$ defined in \cite[pag.\,144 (3)]{Luk69} and separating the poles of $\prod_{i=1}^n\Gm(1-a_i+s)$ from the poles of $\prod_{i=1}^m\Gm(b_i-s)$,
\item $T(m,z)$, with $m\in\N^*$, is the specialization
\begin{align*}T(m,z):&=G_{m-1,m}^{m,0}\left(z\left|\begin{array}{c}
0,\dots, 0\\
-1,-1,\dots,-1,-1
\end{array}\right.\right)\\&=\frac{1}{2\pi\sqrt{-1}}\int_L\left(\frac{-1}{t+1}\right)^{m-1}\Gm(-1-t)z^t{\rm d}t,
\end{align*}where $L$ surrounds the multipole at $t=-1$ and the simple poles at $t=-1+k$ for $k\in\N^*$,
\item $g_m(z)$, with $m\in\N$, is the function
\[g_m(z):=\log^mz\,\Gm(0,z)+mz\sum_{i=0}^{m-1}\frac{(m-1)!}{(m-i-1)!}\log^{m-i-1}z\, T(3+i,z),
\]
\end{itemize}
We have the following explicit formula for $\mc E_k$.
\begin{prop}
For any $j\in\N$, we have
\[\mc E_j(z)=e^z\left(\dl_{0j}-j!\sum_{n=0}^{j-1}\frac{\bt_n g_{j-n-1}(z)}{n!(j-n-1)!}\right).
\]
\end{prop}
\proof
If $f(t)=\sum_{m\geq 1}\frac{a_m}{m!}t^m$ is a formal power series, then $\exp f(t)=\sum_{k=0}^\infty B_k(a_1,\dots,a_k)\frac{t^k}{k!}$, see \cite[pag.\,134, Sec.\,3.3]{Com74}. Consequently, from the Legendre series $\log\Gm(1-t)=\cem t+\sum_{n=2}^\infty\frac{\zeta(n)}{n}t^n$, we deduce
$\frac{1}{\Gm(t)}=\sum_{n=0}^\infty\frac{\bt_n}{n!}t^{n+1}.
$ By results of \cite[\S 4]{GGMS90}, we also have $\der_s^k\Gm(s,z)|_{s=0}=g_k(z)$ for $k\geq 0$. Hence, we have
\[\mc E_j(z)=\left.\der_s^j\left[e^z\left(1-\frac{\Gm(s,z)}{\Gm(s)}\right)\right]\right|_{s=0}=\left.\der_s^j\left[e^z\left(1-\sum_{n=0}^\infty\sum_{k=0}^\infty\frac{\bt_ng_k(z)}{n!k!}s^{n+k+1}\right)\right]\right|_{s=0},
\]and the result follows.
\endproof
\begin{example}The first elements of the sequence $(\mc E_k)_{k\geq 0}$ are 
\begin{align*}
\mc E_0(z)=&\,e^z,&\\
\mc E_1(z)=& -e^z\,\Gm(0,z),&\\
\mc E_2(z)=&-2 e^z \left(G_{2,3}^{3,0}\left(z\left|
\begin{array}{c}
 1,1 \\
 0,0,0 \\
\end{array}
\right.\right)+(\cem+\log z)  \Gamma (0,z)\right),&\\
\mc E_3(z)=&-\frac{1}{2} e^z \Bigg(12 G_{3,4}^{4,0}\left(z\left|
\begin{array}{c}
 1,1,1 \\
 0,0,0,0 \\
\end{array}
\right.\right)+12 (\log z+\cem ) G_{2,3}^{3,0}\left(z\left|
\begin{array}{c}
 1,1 \\
 0,0,0 \\
\end{array}
\right.\right)\\
&+\left(6 \log z (\log z+2 \cem )+6 \cem ^2-\pi ^2\right) \Gamma (0,z)\Bigg),
\end{align*}
\begin{align*}
\mc E_4(z)=&\,2 e^z \Bigg(-12 G_{4,5}^{5,0}\left(z\left|
\begin{array}{c}
 1,1,1,1 \\
 0,0,0,0,0 \\
\end{array}
\right.\right)+\left(-6 \log z (\log z+2 \cem )-6 \cem ^2+\pi ^2\right)\\
&\cdot G_{2,3}^{3,0}\left(z\left|
\begin{array}{c}
 1,1 \\
 0,0,0 \\
\end{array}
\right.\right)-12 (\log z+\cem ) G_{3,4}^{4,0}\left(z\left|
\begin{array}{c}
 1,1,1 \\
 0,0,0,0 \\
\end{array}
\right.\right)+\\
&\Gamma (0,z) \left(-(\log z+\cem ) \left(2 \log z (\log z+2 \cem )+2 \cem ^2-\pi ^2\right)-4 \zeta (3)\right)\Bigg).
\end{align*}
\end{example}

\subsection{Main theorems} 

Let $X_1,\dots, X_h$, with $h\geq 1$, be Fano smooth projective complex varieties. Assume that $\det TX_j=L_j^{\otimes \ell_j}$, with $\ell_j\in\N^*$, for ample line bundles $L_j\to X$.
\vskip2mm
Set $X:=\prod_{j=1}^hX_j$, and $P:=\Pb\left(\mc O_X\oplus\boxtimes_{j=1}^hL_j^{\otimes(-d_j)}\right)$, with $0<d_j<\ell_j$.
\vskip2mm
Let $\dl_j\in H^2(X_j,\C)$, with $j=1,\dots, h$, and let $$\bm\dl=\sum_{j=1}^h1\otimes\dots\otimes\dl_j\otimes\dots\otimes 1\in H^2(X,\C)$$ be the corresponding point of $H^2(X,\C)$, under K\"unneth isomorphism. Denote by $\mc S_{\dl_j}(X_j)$ and $\mc S_{\pi^*\bm\dl}(P)$ the corresponding space of master functions.
\begin{thm}\label{mt2}
The space $\mc S_{\pi^*\bm\dl}(P)$ is contained in the finite dimensional $\C$-vector space generated by the images of the maps $\mathscr B_{\bm\al,\bm\bt,k}\colon\bigotimes_{j=1}^h\mc S_{\dl_j}(X_j)\to\mc O\left(\widetilde{\C^*}\right)$ defined by
\[\Phi_1\otimes\dots\otimes\Phi_h\mapsto \mathscr B_{\bm\al,\bm\bt}\left[\Phi_1,\dots,\Phi_h,\mc E_k\right],
\]where
\[\bm\al=\left(\frac{\ell_1^2}{d_1(d_1-\ell_1)},\dots,\frac{\ell_h^2}{d_h(d_h-\ell_h)},\frac{1}{2}\right),
\qquad \bm\bt=\left(-\frac{d_1}{\ell_1},\dots,-\frac{d_h}{\ell_h},1\right),
\]
and \[k=0,\dots, \dim_\C X+1.
\]In other words, every element of $\mc S_{\pi^*\bm\dl}(P)$ is a finite sum of integrals of the form
\[\frac{1}{2\pi\sqrt{-1}}\int_{\frak H}\prod_{j=1}^h\Phi_j\left(z^{\frac{\ell_j-d_j}{\ell_j}}\la^\frac{d_j}{\ell_j}\right)\mc E_k\left(z^2\la^{-1}\right)e^\la\frac{{\rm d}\la}{\la},\qquad k=0,\dots,\dim_\C X+1,
\]and $\Phi_j\in\mc S_{\dl_j}(X_j)$.
\end{thm}

We can make the result more explicit, in the case of projective bundles over product of projective spaces  \cite{QR98,CMR00,AM00,Str15}.
\begin{thm}\label{mt2.2}
Let $h\geq 1$, and $P$ be the projective bundle $\Pb\left(\mc O_X\oplus \mc O_X(-d_1,\dots,-d_h)\right)$ over the variety $X=\Pb^{n_1-1}\times\dots\times\Pb^{n_h-1}$, with $0<d_i<n_i$ for any $i=1,\dots, h$. Any master function in $\mc S_0(P)$ is a $\C$-linear combination of integrals of the form
\[H_{\bm j}(z):=\frac{1}{(2\pi\sqrt{-1})^{h+1}}\int_{\times \gm_i}\int_\frak H\left[\prod_{i=1}^h\Gm(s_i)^{n_i}\phi^i_{j_i}(s_i)\right]\mc E_k\left(\frac{z^2}{\la}\right)\frac{z^{\sum_i(d_i-n_i)s_i}e^\la}{\la^{1+\sum_id_is_i}}\,{\rm d}\la\,{\rm d}s_1\dots{\rm d}s_h,
\]for $k=0,1,\dots,1-h+\sum_{i=1}^hn_i$, and $\bm j=(j_1,\dots, j_h)\in\prod_{i=1}^h\left\{0,\dots, n_i-1\right\}$. The paths $\gm_i$, with $i=1,\dots,h$, are parabolas of the form ${\rm Re}\,s_i=-\rho_{1,i}({\rm Im}\,s_i)^2+\rho_{2,i}$, for suitable $\rho_{1,i},\rho_{2,i}\in\R_+$, so that they encircle the poles of the factors $\Gm(s_i)^{n_i}$. The function $\phi^i_{j_i}$ is defined as follows
\begin{itemize}
\item for $n_i$ even:
\[\phi^i_{j_i}(s_i):=\exp(2\pi\sqrt{-1}j_is_i),\quad j_i=0,\dots,n_i-1;
\]
\item for $n_i$ odd:
\[\phi^i_{j_i}(s_i):=\exp(2\pi\sqrt{-1}j_is_i+\pi\sqrt{-1}s_i),\quad j_i=0,\dots,n_i-1.
\]
\end{itemize}
\end{thm}
\proof
The result follows by applying Theorem \ref{mt2} to the case $X_i=\Pb^{n_i-1}$, $\ell_i=n_i$. For each factor $\Pb^{n_i-1}$ the space $\mc S_0(\Pb^{n_i-1})$ of master functions equals the space of solutions $g$ of the differential equation
\[\vartheta^{n_i}g=(n_iz)^{n_i}g,\qquad \vartheta=z\frac{d}{dz}.
\]So, a basis of the space $\mc S_0(\Pb^{n_i-1})$ is given by the integrals 
\[g^i_{j_i}(z):=\frac{1}{2\pi\sqrt{-1}}\int_{\gm_i}\Gm(s_i)^{n_i}z^{-n_is_i}\phi^i_{j_i}(s_i){\rm d}s_i,\quad j_i=0,\dots, n_i-1.
\]See e.g.\,\,\cite[Lemma 5]{Guz99}, \cite[Ch.\,6]{CDG18}. For
\[\bm\al=\left(\frac{n_1^2}{d_1(d_1-n_1)},\dots,\frac{n_h^2}{d_h(d_h-n_h)},\frac{1}{2}\right),\quad\bm\bt=\left(-\frac{d_1}{n_1},\dots,-\frac{d_h}{n_h},1\right),
\]we have
\[\mathscr B_{\bm\al,\bm\bt}\left[g^1_{j_1},\dots,g^h_{j_h},\mc E_k\right]=H_{\bm j}(z),
\]where $k=0,\dots, \dim_\C X+1$, and $j_i=0,\dots,n_i-1$ with $i=1,\dots,h$.
\endproof

\section{Proof of the main theorem}\label{sec3}
\subsection{Topological--enumerative solution, and $J$-function}
\begin{defn}Define the functions $\theta_{\beta,p}(z,\bm t),\,\theta_{\beta}(z,\bm t)$, with $\beta=0,\dots, n$ and $p\in\mathbb N$, by
\begin{gather}
\theta_{\beta,p}(\bm t):=\left.\frac{\partial^2\mathcal F_0^X(\bm t^\bullet)}{\partial t^{0,0}\partial t^{\beta,p}}\right|_{\substack{t^{\alpha, p}=0\text{ for }p>1,\\ t^{\alpha,0}=t^\alpha\text{ for }\alpha=0,\dots, n}},\\
\theta_\beta(z,\bm t):=\sum_{p=0}^\infty\theta_{\beta,p}(\bm t)z^p.
\end{gather}
Define the matrix $\Theta(z,\bm t)$ by
\beq
\Theta(z,\bm t)^\alpha_\beta:=\eta^{\alpha\lambda}\frac{\partial\theta_\beta(z,\bm t)}{\partial t^\lambda},\quad \alpha,\beta=0,\dots, n.
\eeq
Let $R$ be the matrix associated with the $\C$-linear morphism $$H^\bullet(X,\C)\to H^\bullet(X,\C),\quad v\mapsto c_1(X)\cup v,$$
with respect to the basis $(T_0,\dots, T_n)$, and define the matrix $Z_{\rm top}(z,\bm t)$ by
\beq
Z_{\rm top}(z,\bm t)\Theta(z,\bm t)z^\mu z^R.
\eeq
\end{defn}

\begin{thm}[\cite{Dub99,CDG20}]
The matrix $Z_{\rm top}(z,\bm t)$ is a fundamental system of solutions of the joint system \eqref{eq1}-\eqref{qde}. \qed
\end{thm}

\begin{defn}
The solution $Z_{\rm top}(z,\bm t)$ is called \emph{topological-enumerative solution} of the joint system \eqref{eq1}, \eqref{qde}.
\end{defn}

\vskip2mm
Let $\hbar$ be an indeterminate.

\begin{defn}
The {\it $J$-function} of $X$ is the $H^\bullet(X,\La_{X,\om})[\![\hbar^{-1}]\!]$--valued function of $\bm\tau\in H^\bullet(X,\C)$ defined by
\[J_X(\bm\tau):=1+\sum_{\al,\la=0}^n\sum_{p=0}^\infty\hbar^{-(p+1)}T_\la\eta^{\al\la}\left.\frac{\der^2\mc F^X_0}{\der t^{0,0}\der t^{\al,p}}\right|_{\substack{t^{\alpha, p}=0\text{ for }p>1,\\ t^{\alpha,0}=\tau^\alpha\text{ for }\alpha=0,\dots, n}}.
\]
\end{defn}

The restriction of the $J$-function to the small quantum locus, i.e.\,\,to points $\bm\tau\in H^2(X,\C)$, has a simpler expansion.
\begin{lem}\cite{CK99,Cot22}\label{sjf}
For $\dl\in H^2(X,\C)$, we have
\[
\pushQED{\qed}
J_X(\dl)=e^\frac{\dl}{\hbar}\left(1+\sum_{\al=0}^n\sum_{\bt\neq 0}\sum_{p=0}^\infty e^{\int_\bt\dl}\langle\tau_pT_\al,1\rangle^X_{0,2,\bt}T^\al\hbar^{-(p+1)}{\bf Q}^\bt\right). \qedhere
\popQED
\]
\end{lem}

\begin{thm}\cite[Th.\,5.2, Cor.\,5.3]{Cot22}
Let $\dl\in H^2(X,\C)$. For $\al=0,\dots, n$, the $(0,\al)$--entry of the matrix $\eta Z_{\rm top}(z,\dl)$ equals
\[z^{\frac{\dim X}{2}}\int_XT_\al\cup J_X(\dl+\log z\cdot c_1(X))\rqh.
\]In particular, the components of the function
\[J_X(\dl+\log z\cdot c_1(X))\rqh,
\]with respect to any $\C$-basis of $H^\bullet(X,\C)$, span the space of master functions $\mc S_\dl(X)$. \qed
\end{thm}

\subsection{The Elezi--Brown theorem} Consider a split vector bundle $V=\bigoplus_{j=0}^rL_j$ over $X$, with line bundles $L_j\to X$ such that
\begin{enumerate}
\item $L_0=\mc O_X$ (not restrictive by Remark \ref{rklines});
\item $L_j^*$ is ample for $j=1,\dots,r$;
\item the class $c_1(X)+c_1(V)$ is ample.
\end{enumerate}

\vskip2mm
Set $P:=\Pb(V)$. Let $s_i\colon X\to P$, with $i=0,\dots, r$, be the section of $\pi\colon P\to X$ determined by the $i$-th summand of $V$. The following lemma describes the Mori cone of $P$.

\begin{lem}[{\cite[Lem.\,\,1.0.1]{Ele05}}]\label{lemc}Assume that each line bundle $L_i^*$, $i=1,\dots, r$, is nef.
\begin{enumerate}
\item If $(T_1,\dots, T_k)$ is a nef basis of $H^2(X,\Q)$, then $(T_1,\dots, T_k,\xi)$ is a nef basis of $H^2(P,\Q)$.
\item The Mori cone of $X$ and $P$ are related by
\[{\rm NE}(P)={\rm NE}(X)\oplus \Z_{\geq 0}\cdot [\ell],
\]where $[\ell]$ is the class of a line in the fiber of $\pi$. Here ${\rm NE}(X)$ is embedded in ${\rm NE}(P)$ via the section $s_{0}$.\qed
\end{enumerate}
\end{lem}

Consider the small $J$-function of $X$: by Lemma \ref{sjf}, it is of the form
\[J_X(\dl)=e^\frac{\dl}{\hbar}\sum_{\bt\in H_2(X,\Z)}J^X_\bt(\dl){\bf Q}^\bt,
\]for suitable coefficients $J_\bt^X(\dl)$.

For each $\bt\in H_2(X,\Z)$ and each $\nu\in\N$, introduce the {\it twisting factor}
\[\mc T_{\nu,\bt}:=\prod_{j=0}^r\frac{\prod_{m=-\infty}^0\left(\xi+\pi^*c_1(L_j)+m\hbar\right)}{\prod_{m=-\infty}^{\nu+\langle c_1(L_j),\bt \rangle}\left(\xi+\pi^*c_1(L_j)+m\hbar\right)},
\]and define the {\it small $I$-function} (defined on $H^2(P,\C)\cong H^2(X,\C)\oplus\,\C\xi$) as the hypergeometric modification
\[I_P(\pi^*\dl+t\xi):=\exp\left(\frac{\pi^*\dl+t\xi}{\hbar}\right)\sum_{\bt}\sum_{\nu\geq 0}\pi^*J^X_\bt(\dl)\cup \mc T_{\nu,\bt}\,{\bf Q}^{\bt+\nu},\quad t\in\C,
\]where the sum $\bt+\nu$ denotes the element $\bt+\nu[\ell]$ of ${\rm NE}(P)$, in the notations of Lemma \ref{lemc}.
\vskip2mm
The following result was first conjectured, and proved in several cases, by A.\,Elezi \cite{Ele05,Ele07}, and finally proved in full generality by J.\,Brown \cite{Bro14}.

\begin{thm}\label{EBT}
If the line bundles $L_j$ satisfy the assumption (1),(2),(3) above, then 
\beq\label{i=j}
\pushQED{\qed}
J_P(\pi^*\dl+t\xi)=I_P(\pi^*\dl+t\xi).\qedhere
\popQED
\eeq
\end{thm}
\begin{rem}
The result proved by J.\,Brown in \cite{Bro14} is actually more general, since it concerns the more general case of a toric fiber bundle on $X$ rather than a projective bundle. Its formulation, however, is less explicit than the equality \eqref{i=j}. Namely, Brown's result is stated only in terms of 
the A.\,Givental's formalism of Lagrangian cones \cite{Giv04,CG07}.
\end{rem}

\subsection{The Ribenboim's algebras \texorpdfstring{$\mathscr F_{\bm\kappa}(A)$}{}}Let $h\in\N^*$, and fix a tuple $\bm \kappa:=(\kappa_1,\dots,\kappa_h)\in(\C^*)^h$. In this and in the next section, let $(A,+,\cdot,1_A)$ be an associative, commutative, unital and finite dimensional $\C$-algebra, with nilradical ${\rm Nil}(A):=\{a\in A\colon\exists \,n\in\N\,\,\text{s.t.\,\,}a^n=0\}$.
\vskip2mm
Set $\N_A:=\{n\cdot 1_A\colon n\in\N\}$. Define the monoid $M_{A,{\bm\kappa}}$ as the (external) direct sum of monoids
\[M_{A,{\bm\kappa}}:=\left(\bigoplus_{j=1}^h\kappa_j\mathbb N_A\right)\oplus{\rm Nil}(A).\]
If $x\in M_{A,\bm\kappa}$, denote by $x'$ its ``nilpotent part'', i.e. the projection of $x$ onto ${\rm Nil}(A)$. 
\vskip2mm
We have two maps $\nu_{\bm\kappa}\colon M_{A,{\bm\kappa}}\to \mathbb N^h$ and $\iota_{\bm\kappa}\colon M_{A,{\bm\kappa}}\to A$ defined by
\[\nu_{\bm\kappa}((\kappa_in_i1_A)_{i=1}^h,r):=(n_i)_{i=1}^h,\quad \iota_{\bm\kappa}((\kappa_in_i1_A)_{i=1}^h,r):=\sum_{i=1}^h\kappa_in_i1_A+r.
\]By universal property of the direct sums of monoids, the natural inclusions $M_{A,\kappa_i}\to M_{A,{\bm\kappa}}$ induce a unique morphism
\[\rho_{\bm\kappa}\colon \bigoplus_{i=1}^hM_{A,\kappa_i}\to M_{A,\bm\kappa}.
\]

On $M_{A,{\bm\kappa}}$ we can define the partial order 
\[x\leq y \quad\text{iff}\quad x'=y'\text{ and }\nu_{\bm\kappa}(x)\leq \nu_{\bm\kappa}(y),
\]the order on $\mathbb N^h$ being the lexicographical one. This order makes $(M_{A,{\bm\kappa}},\leq)$ a strictly ordered monoid, that is 
\[\text{if } a,b\in M_{A,\bm\kappa}\text{ are such that }a<b,\,\,\text{then } a+c<b+c\quad\text{for all } c\in M_{A,\bm\kappa}.
\]
Define $\mathscr F_{\bm\kappa}(A)$ to be the set of all functions $f\colon M_{A,{\bm\kappa}}\to A$ whose support
\[{\rm supp}(f):=\left\{a\in M_{A,{\bm\kappa}}\colon f(a)\neq 0\right\}
\]is 
\begin{enumerate}
\item {\it Artinian}, i.e.\,\,every subset of ${\rm supp}(f)$ admits a minimal element,
\item and {\it narrow}, i.e.\,\,every subset of ${\rm supp}(f)$ of pairwise incomparable elements is finite.
\end{enumerate}
The set $\mathscr F_{\bm \kappa}(A)$ is an $A$-module with respect to pointwise addition and multiplication of $A$-scalars.
We will denote the element $f\in\mathscr F_{\bm \kappa}(A)$ by
\[f=\sum_{a\in M_{A,\bm\kappa}}f(a)Z^a,
\]where $Z$ is an indeterminate. Given $f_1,f_2\in\mathscr F_{\bm \kappa}(A)$, define
\[f_1\cdot f_2:=\sum_{s\in M_{A,\bm\kappa}}\left(\sum_{(p,q)\in X_s(f,g)}f_1(p)\cdot f_2(q)\right)Z^{s},
\]where we set
\[X_s(f,g):=\left\{(p,q)\in M_{A,\bm\kappa}\times M_{A,\bm\kappa}\colon p+q=s,\quad f_1(p)\neq 0,\quad f_2(q)\neq 0\right\}.
\]

The following result is a consequence of P.\,Ribenboim's theory of generalized power series \cite{Rib92,Rib94}.
\begin{thm}
The product above is well-defined. The set $\mathscr F_{\bm \kappa}(A)$ is equipped with an $A$-algebra structure with respect to the operations above.\qed
\end{thm}

\begin{defn}
Let $r_o\in{\rm Nil}(A)$. We say that an element $f\in \mathscr F_{\bm \kappa}(A)$ is \emph{concentrated at $r_o$} if 
\[{\rm  supp}(f)\subseteq \left(\bigoplus_{i=1}^h\kappa_i\N_A\right)\times\{r_o\}.
\]
\end{defn}

\subsection{Formal Borel multitransform} Given two $h$-tuples $\bm \al,\bm\bt\in(\C^*)^h$, we set $\bm\al\cdot\bm\bt:=(a_i\bt_i)_{i=1}^h$, and $\bm\al^{-1}:=\left(\frac{1}{\al_i}\right)_{i=1}^h$.

\begin{defn}\label{fs}
Let $F\in\C[\![x]\!]$ be a formal power series $F(x)=\sum_{k=0}^\infty a_kx^k$. For $\al\in{{\rm Nil}(A)}$ define $F(\al)\in A$ by the finite sum 
\[F(\al)=\sum_{k=0}^\infty a_k\al^k.
\]
If $F$ is invertible, i.e. $a_0\neq 0$, then $F(\al)$ is invertible in $A$.
\end{defn}
In what follows we will usually take $F(x)=\Gamma(\lambda+x)$ with $\lambda\in\C\setminus\mathbb Z_{\leq 0}$, or $F(x)=1/\Gm(\la+x)$ with\footnote{Recall that the function $1/\Gm(s)$ is an entire function.} $\la\in\C$, where $\Gamma$ denotes the Euler Gamma function.

\begin{defn}
Let ${\bm\alpha,\bm \beta, \bm \kappa}\in(\mathbb C^*)^h$. We define the \emph{Borel $(\bm\alpha,\bm\beta)$-multitransform} as the $A$-linear morphism
\[\mathscr B_{\bm\alpha,\bm\beta}\colon \bigotimes_{j=1}^h\mathscr F_{\kappa_j}(A)\to \mathscr F_{\bm\al^{-1}\cdot\bm\bt^{-1}\cdot\bm\kappa}(A),
\]which is defined, on decomposable elements, by
\[\mathscr B_{\bm\alpha,\bm\beta}\left(\bigotimes_{j=1}^h\left(\sum_{s_j\in M_{A,\kappa_j}}f_{s_j}^jZ^{s_j}\right)\right):=\sum_{\substack{s_j\in M_{A,\kappa_j}\\ j=1,\dots, h}}\frac{\prod_{i=1}^hf^i_{s_i}}{\Gamma\left(1+\sum_{\ell=1}^h\iota_{ \kappa_\ell}(s_\ell)\beta_\ell\right)}Z^{\rho_{\bm \kappa}\left(\oplus_{\ell=1}^h \frac{s_\ell}{\al_\ell\bt_\ell }\right)}.
\]
\end{defn}

Let $s=((\kappa_in_i1_A)_{i=1}^h,r)\in M_{A,\bm\kappa}$. We define the {\it analytification} $\widehat{Z^s}$ of the monomial $Z^s\in\mathscr F_{\bm \kappa}(A)$ to be the $A$-valued holomorphic function 
\[\widehat{Z^s}\colon \widetilde{\C^*}\to A,\quad \widehat{Z^s}(z):=z^{\sum_{i=1}^h\kappa_i n_i}\sum_{j=1}^\infty\frac{r^j}{j!}\log^j z.
\]Notice that the sum is finite, since $r\in{\rm Nil}(A)$.

Let $f\in\mathscr F_{\bm \kappa}(A)$ be a series
\[f(Z)=\sum_{s\in M_{A,\bm\kappa}}f_aZ^s,\quad\text{such that}\quad  {\rm card\ supp}(f)\leq \aleph_0.
\]

The {\it analytification} $\hat f$ of $f$ is the $A$-valued holomorphic function defined, if the series absolutely converges, by
\[\widehat f\colon W\subseteq \widetilde{\C^*}\to A,\quad \widehat f(z):=\sum_{s\in M_{A,\bm\kappa}}f_a\widehat{Z^s}(z).
\]

\begin{thm}\label{BLAF}\cite[Th.\,6.13]{Cot22}
Let $f_i\in \mathscr F_{ \kappa_i}(A)$ such that 
\begin{itemize}
\item ${\rm card\ supp}(f_i)\leq\aleph_0$ for $i=1,\dots, h$,
\item the functions $\widehat f_i$ are well defined on $\C\setminus\R_{<0}$.
\end{itemize}
We have
\[\reallywidehat{\mathscr B_{\bm\al,\bm\bt}[\bigotimes_{j=1}^hf_j]}=\mathscr B_{\bm\al,\bm\bt}[\widehat f_1,\dots, \widehat f_h],
\]
provided that both sides are well-defined.\qed
\end{thm}

In the remaining part of the paper, we will consider the Ribenboim algebras $\mathscr F_{\bm\kappa}(X):=\mathscr F_{\bm\kappa}(H^\bullet(X,\C))$ where $X$ is a smooth projective variety such that $H^{\rm odd}(X,\C)=0$. This is necessary in order to work with commutative cohomological algebras.
\subsection{The Ribenboim $\mc E$-series of a projective bundle}
Given a projective bundle $P\to X$ as in Section \ref{secpb}, we associate to $P$ a distinguished element of $\mathscr F_1(P)$ as follows.
\begin{defn}
The $\mc E$-series of $P$ is the formal power series $\mc E_P(\xi;Z)\in \mathscr F_1(P)$ defined by
\[\mc E_P(\xi;Z):=\sum_{k=0}^\infty \frac{Z^{k+\xi}}{\Gm(1+k+\xi)}.
\] 
\end{defn}

\begin{prop}\label{ecomp}
The analytification $\widehat{\mc E_P}(\xi;z)$ equals
\[\widehat{\mc E_P}(\xi;z)=\sum_{k\geq 0}\frac{1}{k!}\mc E_k(z)\xi^k.
\]
\end{prop}
\proof
It follows from Proposition \ref{expe}, the definitions of the $\mc E_k$-functions, Definition \ref{fs} and that of analytification of a Ribenboim series.
\endproof

\subsection{Proof of Theorem \ref{mt2}}
Set $\rho_j:=c_1(L_j)$, with $j=1,\dots,h$.

By K\"unneth isomorphism, and by the universal property of coproduct of algebras (i.e. tensor product), we have injective\footnote{In particular, we have inclusions $\mathscr F_{\bm k}(X_j)\to \mathscr F_{\bm k}(X)$. } maps $H^\bullet(X_i,\C)\to H^\bullet(X,\C)$. In order to ease the computations, in the next formulas we will not distinguish an element of $H^\bullet(X_i,\C)$ with its image in $H^\bullet(X,\C)$. The same will be applied to elements of $H_2(X,\Z)$. So, for example, we have
\[c_1(X)=\sum_{j=1}^h\ell_j\rho_j,\qquad c_1(P)=\pi^*\left(\sum_{j=1}^h(\ell_j-d_j)\rho_j\right)+2\xi,
\]where the last identity follows from Lemma \ref{lccp}.

The space of master functions $\mc S_{\bm\dl}(X)$ is generated by the components (with respect an arbitrary basis of $H^\bullet(X,\C)$) of the small $J$-function
\[J_{X}\left(\bm\delta+c_1(X)\log z\right)\rqh=\bigotimes_{j=1}^h J_{X_j}(\delta_j+c_1(X_j)\log z)\rqh,
\]by the R.\,Kaufmann's quantum K\"unneth formula \cite{Kau96}.
By Lemma \ref{sjf}, for each $j=1,\dots,h$, we have
\begin{multline*}
\left.J_{X_j}(\delta_j+\log z\cdot c_1(X_j))\right\rqh\\=e^{\delta_j} z^{c_1(X_j)}\left(1+\sum_\al\sum_{\substack{\bt\in{\rm NE}(X_j)\\\bt\neq 0}}\sum_{k=0}^\infty e^{\int_\bt\delta_j} z^{\int_\bt c_1(X_j)}\langle\tau_k T_{\al,j},1\rangle_{0,2,\bt}^{X_j} T_j^\al\right),
\end{multline*}
where $(T_{0,j},\dots, T_{n_j,j})$ is a fixed basis of $H^\bullet(X_j,\C)$, $T^\al_j:=\sum_{\la=0}^{n_j}\eta^{\al\la}_jT_{\al,j}$, and $\eta_j$ is the Poincar\'e metric on $H^\bullet(X_j,\C)$. So, we can rewrite the $J$-function $J_{X}\left(\bm\delta+c_1(X)\log z\right)\rqh$ in the form
\beq\label{jx2}
J_{X}\left(\bm\delta+c_1(X)\log z\right)\rqh=\bigotimes_{j=1}^h\sum_{k_j=0}^\infty J^j_{k_j}(\dl_j)z^{k_j\ell_j+\ell_j\rho_j},
\eeq
where the coefficient $J^j_{k_j}(\dl_j)$ equals
\[J^j_{k_j}(\dl_j)=e^{\dl_j}\sum_{\al=0}^{n_j}\sum_{p=0}^\infty\sum_{\substack{\bt\in{\rm NE}(X_j)\\ \langle\rho_j,\bt\rangle=k_j}}e^{\int_\bt\dl_j}\langle\tau_pT_{\al,j},1\rangle^{X_j}_{0,2,\bt}T_j^\al.
\]Analogously, the space of master functions $\mc S_{\pi^*\bm\dl}(P)$ is spanned by the components of the $H^\bullet(P,\C)$-valued function
$J_P(\pi^*\bm\dl+c_1(P)\log z)\rqh$. 
\vskip2mm
By Elezi--Brown Theorem \ref{EBT}, we have
\begin{multline*}
J_P(\pi^*\bm \dl+c_1(P)\log z)\rqh=I_P\left[\pi^*\left(\bm\dl+\log z(c_1(X)+c_1(\boxtimes_{j=1}^hL_j^{\otimes(-d_j)})\right)+2\xi \log z\right]\rqh\\
=\sum_{k_1,\dots,k_h\in\N}\pi^*\left(\bigotimes_{j=1}^h J^j_{k_j}\left(\dl_j+\log z(\ell_j-d_j)\rho_j\right) \right)\frac{z^{2\nu+2\xi}}{\prod_{m=1}^{\nu}(\xi+m)\prod_{m=1}^{\nu-\sum_jk_jd_j}(\xi-\sum_jd_j\rho_j+m)}\\
=\sum_{k_1,\dots,k_h\in\N}\pi^*\left(\bigotimes_{j=1}^h J^j_{k_j}\left(\dl_j\right) \right)\frac{\Gm(1+\xi)}{\Gm(1+\nu+\xi)}
\frac{\Gm(1+\xi-\sum_jd_j\rho_j)}{\Gm(1+\nu+\xi-\sum_jd_j(k_j+\rho_j))}z^{(\ell_j-d_j)(k_j+\rho_j)+2\nu+2\xi}.
\end{multline*}
On the one hand, each factor in the tensor product \eqref{jx2} equals the analytification $\widehat{{\rm J}_{X_j}}$ of the series ${\rm J}_{X_j}\in\mathscr F_{\ell_j}(X)$ defined by
\[{\rm J}_{X_j}(Z)=\sum_{k_j=0}^\infty J^j_{k_j}(\dl_j)Z^{k_j\ell_j+c_1(X_j)}.
\]On the other hand, consider the Borel multitransform $\mathscr B_{\bm\al,\bm\bt}\left[\pi^*{\rm J}_{X_1},\dots,\pi^*{\rm J}_{X_h},\mc E_P\right]$ for arbitrary tuples $\bm\al=(\al_1,\dots,\al_{h+1})$ and $\bm\bt=(\bt_1,\dots,\bt_{h+1})$. We have
\begin{multline*}
\mathscr B_{\bm\al,\bm\bt}\left[\pi^*{\rm J}_{X_1},\dots,\pi^*{\rm J}_{X_h},\mc E_P\right]\\=
\sum_{k_1,\dots,k_h=0}^\infty\sum_{k=0}^\infty\frac{\pi^*\left(\bigotimes_{j=1}^h J^j_{k_j}\left(\dl_j\right) \right)}{\Gm(1+k+\xi)}\frac{Z^{\sum_{j=1}^h\frac{k_j\ell_j+c_1(X_j)}{\al_j\bt_j}+\frac{k+\xi}{\al_{h+1}\bt_{h+1}}}}{\Gm(1+\sum_{j=1}^h\bt_j(k_j\ell_j+\ell_j\rho_j)+\bt_{h+1}(k+\xi))},
\end{multline*}
by definitions of the formal Borel multitransform and of the $\mc E$-series. Then, for the choice of weights
\[\al_j=\frac{\ell_j^2}{d_j(d_j-\ell_j)},\quad j=1,\dots, h,\qquad \al_{h+1}=\frac{1}{2},
\]
\[\bt_j=-\frac{d_j}{\ell_j},\quad j=1,\dots, h,\qquad \bt_{h+1}=1,
\]we have
\[J_P(\pi^*\bm \dl+c_1(P)\log z)\rqh=\Gm(1+\xi)\Gm\left(1+\xi-\sum_jd_j\rho_j\right)\reallywidehat{\mathscr B_{\bm\al,\bm\bt}\left[\pi^*{\rm J}_{X_1},\dots,\pi^*{\rm J}_{X_h},\mc E_P\right]}.
\]
Then, the space $\mc S_{\pi^*\bm\dl}(P)$ is spanned by the components, with respect to an arbitrary basis of $H^\bullet(P,\C)$, of the analytification $\reallywidehat{\mathscr B_{\bm\al,\bm\bt}\left[\pi^*{\rm J}_{X_1},\dots,\pi^*{\rm J}_{X_h},\mc E_P\right]}$. This follows from the invertibility of  the morphism
\[H^\bullet(P,\C)\to H^\bullet(P,\C),\quad v\mapsto \Gm(1+\xi)\Gm\left(1+\xi-\sum_jd_j\rho_j\right)v.
\]The statement of Theorem \ref{mt2} then follows from Theorem \ref{BLAF} and Proposition \ref{ecomp}.

\section{An example: blowing-up a point in $\Pb^2$}\label{sec4}
\subsection{Classical and quantum cohomology} Let $P$ be the blow-up of a point in $\Pb^2$. The variety $P$ admits a $\Pb^1$-bundle structure on $\Pb^1$, namely 
\[P=\Pb(\mc O\oplus\mc O(-1)).
\]By looking at $P$ as the compactification of the total space of $\mc O(-1)\to\Pb^1$, we can introduce two cycles $\varepsilon,\nu\in H_2(P,\Z)$ defined as follows:
\begin{itemize}
\item $\varepsilon$ represents the homology class of the section ``at infinity'' of $\mc O(-1)$,
\item $\nu$ represents the homology class of a fiber of $\mc O(-1)$.
\end{itemize}The classical cohomology $H^\bullet(P,\C)$ admits the $\C$-basis 
\[T_0:=1,\quad T_1:={\rm PD}(\varepsilon),\quad T_2:={\rm PD}(\nu),\quad T_3:={\rm PD}(\rm pt),
\]where ${\rm PD}(\al)$ denotes the Poincar\'e dual class of $\al\in H_\bullet(X,\C)$. We denote by $(t^0,t^1,t^2,t^3)$ the corresponding dual coordinates.

The classical cohomology algebra admits the following presentation
\[H^\bullet(P,\C)\cong \frac{\C[T_1,T_2]}{\langle T_2^2,\,T_1^2-T_1T_2\rangle},\quad \text{and }T_1T_2=T_3.
\]The small quantum cohomology algebra, at a point $\tau=t^1T_1+t^2T_2\in H^2(P,\C)$, admits the following presentation
\[QH^\bullet_\tau(P)\cong \frac{\C[T_1,T_2,q_1,q_2]}{\langle T_2^2-(T_1-T_2)q_2,\,T_1^2-T_1T_2-q_1\rangle},\quad q_1=\exp t^1,\,q_2=\exp t^2.
\]See, for example, \cite{Aud97}.
\subsection{QDE and $\La$-matrix}The quantum differential equation \eqref{qde} of $P$ at a point $\tau=t^1T_1+t^2T_2\in H^2(P,\C)$ is 
\beq\label{qdeP}\frac{d}{dz}\sig(\tau,z)=\left(\mc U(\tau)+\frac{1}{z}\mu(\tau)\right)\sig(\tau,z),
\eeq
where $\mu$ is the grading operator on $H^\bullet(P,\C)$, and $\mc U$ denotes the operator of $\circ_\tau$-quantum multiplication by the first Chern class $c_1(P)$. In components, these operators are represented by the matrices
\[\mu(\tau)={\rm diag}(-1,0,0,1),\qquad \mc U(\tau)=\begin{pmatrix}
0&2q_1&0&3q_1^2q_2\\
2&q_1q_2&q_1q_2&0\\
-1&-2q_1q_2&-2q_1q_2&2q_1\\
0&5&2&0
\end{pmatrix}.
\]
At a point $\tau\in H^2(P,\C)$ the determinant of the $\La$-matrix equals
\[\det \La(\tau,z)=-\frac{z}{(27   q_1^{2 }q_2^2+256 q_1) z-24 q_2 q_1},
\]so that the intersection of the $\mc A_\La$-stratum and small quantum locus $H^2(P,\C)$ is 
\[
\label{Alf1}
\mc A_\La\cap H^2(P,\C)=\left\{\tau\colon 27   q_1^{2 }q_2^2+256 q_1=0\right\}.
\]In particular, $0\notin \mc A_\La$ and the $\La$-matrix at the point $0\in H^\bullet(P,\C)$ equals
\[\La(0,z)=\left(
\begin{array}{cccc}
 1 & \frac{204 z^3-105 z^2+9 z+8}{z^2 (283 z-24)} & \frac{-408 z^3-73 z^2+6 z-16}{z^2 (283 z-24)} & \frac{-218 z^3+16 z^2+35 z-6}{z^2 (283 z-24)} \\
 0 & \frac{169 z^2-9 z-8}{z (283 z-24)} & \frac{-55 z^2-6 z+16}{z (283 z-24)} & \frac{-28 z^2-35 z+6}{z (283 z-24)} \\
 0 & \frac{z}{283 z-24} & -\frac{2 z}{283 z-24} & \frac{35 z-3}{283 z-24} \\
 0 & -\frac{8 z}{283 z-24} & \frac{16 z}{283 z-24} & \frac{3 z}{283 z-24} \\
\end{array}
\right),
\]with determinant 
\[\det \La(0,z)=\frac{z}{24-283z}.
\]The quantum differential equation \eqref{qdeP}, specialized at $\tau=0$, can be put in the form \eqref{qde.2} and \eqref{qde.3}. Consequently, it can be reduced to a single scalar differential equation (the master differential equation) in the unknown function $\varpi_0(0,z)=\varpi_0(z)$.

If we introduce the factorization $\varpi_0(z)=z\Phi(z)$, the function $\Phi\in\mc O(\widetilde{\C^*})$ is a master function in the sense of Definition \ref{defmastfun}. The master differential equation, in terms of $\Phi$, reads
\begin{align}\label{qdiffdisp2bis}
&(283 z-24)\vartheta^4\Phi+\left(283 z^2-590 z+24\right)\vartheta^3\Phi+ \left(-2264 z^2+192 z+3\right)\vartheta^2\Phi\\
\nonumber&-4 z^2 \left(2547 z^2+350 z-104\right)\vartheta\Phi+z^2 \left(-3113 z^3-9924 z^2+1476 z+192\right)\Phi=0,
\end{align}
where $\vartheta:=z\frac{d}{dz}$. The space of solutions of this equation coincides with $\mc S_0(P)$.

\subsection{Solutions as Borel multitransforms} According to Theorem \ref{mt2}, we can integrate the differential equation \eqref{qdiffdisp2bis} starting from the space of master functions $\mc S_0(\Pb^1)$. This is the space of functions $\Ups\in\mc O(\widetilde{\C^*})$ solving the differential equation
\[\vartheta^2\Upsilon=4z^2\Upsilon.
\]See e.g.\,\,\cite{Guz99}, \cite[Ch.\,6]{CDG18}. The functions $\Ups$ 
have the following expansion at $z=0$:
\beq\label{expz0}
\Ups(z)=\sum_{k=0}^\infty(a_k+b_k\log z)\frac{z^{2k}}{(k!)^2},
\eeq
where $a_0,b_0$ are arbitrary complex numbers, and all the remaining coefficients are uniquely determined by the difference equations
\[
a_{k-1}=\frac{b_k}{k}+a_k,\qquad b_{k-1}=b_{k},\quad k\geq 1.
\]In particular, we obtain 
\[a_k=a_0-b_0H_k,\qquad H_k:=\sum_{m=1}^k\frac{1}{m}.
\]
Consider the basis $(\Ups_1,\Ups_2)$ of $\mc S_0(\Pb^1)$ defined by the following choices of $(a_0,b_0)\in\C^2$:
\[\Ups_1(z)\text{ is defined by equation \eqref{expz0} with }(a_0,b_0)=(1,0),
\]
\[\Ups_2(z)\text{ is defined by equation \eqref{expz0} with }(a_0,b_0)=(0,1).
\]
Set
\[\bm\al=\left(-4,\frac{1}{2}\right),\qquad \bm\bt=\left(-\frac{1}{2},1\right).
\]
\begin{thm}\label{lastth}
The space $\mc S_0(P)$ of solutions of the differential equation \eqref{qdiffdisp2bis} admits the following $\C$-basis:
\beq\label{basiss0p}
\mathscr B_{\bm\al,\bm\bt}[\Ups_1,\mc E_0],\qquad \mathscr B_{\bm\al,\bm\bt}[\Ups_2,\mc E_0],\qquad \mathscr B_{\bm\al,\bm\bt}[\Ups_1,\mc E_1],\qquad 4\mathscr B_{\bm\al,\bm\bt}[\Ups_2,\mc E_1]+\mathscr B_{\bm\al,\bm\bt}[\Ups_1,\mc E_2].
\eeq
\end{thm}
\proof Let $\Ups$ be as in \eqref{expz0}.
Let us compute the expansions at $z=0$ of the functions $$\mathscr B_{\bm\al,\bm\bt}[\Ups,\mc E_j](z)=\frac{1}{2\pi\sqrt{-1}}\int_\frak H\Phi\left(z^\frac{1}{2}\la^\frac{1}{2}\right)\mc E_j\left(\frac{z^2}{\la}\right)e^\la\frac{{\rm d}\la}{\la},\quad j=0,1,2,$$
where
\begin{align*}&\mc E_0(z)=e^z,\qquad
\mc E_1(z)= -e^z\,\Gm(0,z),\\
&\mc E_2(z)=-2 e^z \left(G_{2,3}^{3,0}\left(z\left|
\begin{array}{c}
 1,1 \\
 0,0,0 \\
\end{array}
\right.\right)+(\cem+\log z)  \Gamma (0,z)\right).
\end{align*}For shortening the following formulas, set
\[C^k_z:=\left(\frac{1}{\Gm}\right)^{(k)}\!\!\!(z),\quad k\geq 0,\,z\in\C.
\]
For $j=0$, we have
\begin{align*}
\mathscr B_{\bm\al,\bm\bt}[\Ups,\mc E_0](z)&=\frac{1}{2\pi\sqrt{-1}}\int_\frak H\sum_{k=0}^\infty\sum_{h=0}^\infty\left(a_k+\frac{b_k}{2}\log z+\frac{b_k}{2}\log \la\right)\frac{z^{k+2h}}{(k!)^2h!}\frac{e^\la}{\la^{1+h-k}}{\rm d}\la\\
&=\sum_{k=0}^\infty\sum_{h=0}^\infty\left[\left(a_k+\frac{b_k}{2}\log z\right)C^0_{1+h-k}
-\frac{b_k}{2}C^1_{1+h-k}
\right]
\frac{z^{k+2h}}{(k!)^2h!},
\end{align*}
where in the last equality we used the identity
\[C^k_z=\left(\frac{1}{\Gm}\right)^{(k)}\!\!\!(z)=\frac{1}{2\pi\sqrt{-1}}\int_\frak He^\la \frac{(-\log\la)^k}{\la^z}{\rm d}\la,\qquad k\geq 0.
\]
From the expansions
\[\Gm(0,z)=-\cem-\log z+\sum_{j=1}^\infty\frac{(-1)^{j+1}}{j!j}z^j,
\]
\[G_{2,3}^{3,0}\left(z\left|
\begin{array}{c}
 1,1 \\
 0,0,0 \\
\end{array}
\right.\right)=\frac{1}{12} \left(6 \log ^2(z)+12 \cem  \log (z)+6 \cem ^2+\pi ^2\right)+\sum_{\ell=1}^\infty\frac{(-1)^{\ell}}{\ell!\ell^2}z^\ell,
\]we deduce
\begin{align*}
&\mathscr B_{\bm\al,\bm\bt}\left[\Ups,\mc E_1\right](z)=\frac{\cem}{2\pi\sqrt{-1}}\int_\frak H\sum_{k=0}^\infty\sum_{h=0}^\infty\left(a_k+\frac{b_k}{2}\log z+\frac{b_k}{2}\log \la\right)\frac{z^{k+2h}}{(k!)^2h!}\frac{e^\la}{\la^{1+h-k}}{\rm d}\la\\
&+\frac{1}{2\pi\sqrt{-1}}\int_\frak H\sum_{k=0}^\infty\sum_{h=0}^\infty\left(a_k+\frac{b_k}{2}\log z+\frac{b_k}{2}\log \la\right)(2\log z-\log\la)\frac{z^{k+2h}}{(k!)^2h!}\frac{e^\la}{\la^{1+h-k}}{\rm d}\la\\
&+\frac{1}{2\pi\sqrt{-1}}\int_\frak H\sum_{k=0}^\infty\sum_{h=0}^\infty\sum_{j=1}^\infty\left(a_k+\frac{b_k}{2}\log z+\frac{b_k}{2}\log \la\right)\frac{(-1)^j}{j\cdot j!h!(k!)^2}\frac{z^{k+2(j+h)}e^\la}{\la^{1+j+h-k}}{\rm d}\la
\end{align*}
\begin{align*}
&=(\cem+2\log z)\sum_{k=0}^\infty\sum_{h=0}^\infty\left[\left(a_k+\frac{b_k}{2}\log z\right)C^0_{1+h-k}
-\frac{b_k}{2}C^1_{1+h-k}
\right]\frac{z^{k+2h}}{(k!)^2h!}\\
&+\sum_{k=0}^\infty\sum_{h=0}^\infty\left[\left(a_k+\frac{b_k}{2}\log z\right)C^1_{1+h-k}
-\frac{b_k}{2}C^2_{1+h-k}
\right]\frac{z^{k+2h}}{(k!)^2h!}\\
&+\sum_{k=0}^\infty\sum_{h=0}^\infty\sum_{j=1}^\infty\left[\left(a_k+\frac{b_k}{2}\log z\right)C^0_{1+j+h-k}
-\frac{b_k}{2}C^1_{1+j+h-k}
\right]\frac{(-1)^jz^{k+2(j+h)}}{j\cdot j!h!(k!)^2}.
\end{align*}
Similarly, we have
\begin{align*}
&\mathscr B_{\bm\al,\bm\bt}\left[\Ups,\mc E_2\right](z)=\\
&\left(\cem^2-\frac{\pi^2}{6}+4\cem\log z+4\log^2z\right)\sum_{k=0}^\infty\sum_{h=0}^\infty\left[\left(a_k+\frac{b_k}{2}\log z\right)C^0_{1+h-k}
-\frac{b_k}{2}C^1_{1+h-k}
\right]\frac{z^{k+2h}}{(k!)^2h!}\\
&+(2\cem+4\log z)\sum_{k=0}^\infty\sum_{h=0}^\infty\left[\left(a_k+\frac{b_k}{2}\log z\right)C^1_{1+h-k}
-\frac{b_k}{2}C^2_{1+h-k}
\right]\frac{z^{k+2h}}{(k!)^2h!}\\
&+\sum_{k=0}^\infty\sum_{h=0}^\infty\left[\left(a_k+\frac{b_k}{2}\log z\right)C^2_{1+h-k}
-\frac{b_k}{2}C^3_{1+h-k}
\right]\frac{z^{k+2h}}{(k!)^2h!}\\
&+(2\cem+4\log z)\sum_{k=0}^\infty\sum_{h=0}^\infty\sum_{j=1}^\infty\left[\left(a_k+\frac{b_k}{2}\log z\right)C^0_{1+j+h-k}
-\frac{b_k}{2}C^1_{1+j+h-k}
\right]\frac{(-1)^jz^{k+2(j+h)}}{j\cdot j!h!(k!)^2}\\
&+2\sum_{k=0}^\infty\sum_{h=0}^\infty\sum_{j=1}^\infty\left[\left(a_k+\frac{b_k}{2}\log z\right)C^1_{1+j+h-k}
-\frac{b_k}{2}C^2_{1+j+h-k}
\right]\frac{(-1)^jz^{k+2(j+h)}}{j\cdot j!h!(k!)^2}\\
&+2\sum_{k=0}^\infty\sum_{h=0}^\infty\sum_{\ell=1}^\infty\left[\left(a_k+\frac{b_k}{2}\log z\right)C^0_{1+\ell+h-k}
-\frac{b_k}{2}C^1_{1+\ell+h-k}
\right]\frac{(-1)^{1+\ell} z^{k+2(\ell+h)}}{\ell^2\cdot \ell!h!(k!)^2}.
\end{align*}
Hence, after some computations, one can explicitly check that the functions in \eqref{basiss0p} solve the differential equation \eqref{qdiffdisp2bis}. This can be done by expressing the coefficients $C^k_z$ in terms of values of polygamma functions $\psi^{(k)}$, namely
\[C^k_z=\frac{1}{\Gm(z)}B_k\left(-\psi^{(0)}(z),\,-\psi^{(1)}(z),\,-\psi^{(2)}(z),\dots,\,-\psi^{(k-1)}(z)\right),
\]where $B_k$ is the $k$-th Bell polynomial as in \eqref{bellp}, and by invoking the following well-known identities (see \cite[Ch.\,5]{OLBC10}):
\begin{align*}
\psi^{(k)}(z+1)&=\psi^{(k)}(z)+\frac{(-1)^kk!}{z^{k+1}},\quad k\geq 0,\\
\psi^{(0)}(n)&=H_{n-1}-\cem,\quad n\geq 1,\quad \psi(z):=\frac{\Gamma'(z)}{\Gamma(z)}.\qedhere
\end{align*}

The explicit formulas of Theorem \ref{lastth} open the possibility to the study of analytical properties of the solutions of the qDE \eqref{qdeP}, such as their asymptotics, Stokes phenomenon, etc. As shown in \cite{Cot22}, these analytical properties remarkably encode information not only about the topology, but even about the algebraic geometry of the surface $P$: the entries of the connection and Stokes matrices can be related to explicit characteristic classes of $P$ and of objects of exceptional collections in the derived category $\mathcal D^b(P)$, see \cite[Th.\,11.8.3]{Cot22}.

We underline that the basis of solutions of the qDE \eqref{qdeP} considered in \cite{Cot22} was constructed via Laplace $(\al,\bt)$-multitransforms, by realizing $P$ as an hypersurface in $\Pb^1\times\Pb^2$. Notice that the computations involved in the proof of Theorem \ref{lastth} are much faster then those performed in \cite[Ch.\,11, App.\,B]{Cot22}.

\endproof

\end{document}